\documentclass[11pt]{amsart}

\usepackage{amsmath}
\usepackage{amssymb}
\usepackage{enumerate}
\usepackage{euscript}
\usepackage{amscd}
\usepackage[all]{xy}

\newtheorem{thm}{Theorem}[section]
\newtheorem{lem}[thm]{Lemma}
\newtheorem{prop}[thm]{Proposition}
\newtheorem{cor}[thm]{Corollary}
\theoremstyle{definition}
\newtheorem{defn}[thm]{Definition}
\newtheorem{rem}[thm]{Remark}

\newtheorem{asmp}[thm]{Assumption}

\numberwithin{equation}{section}
\theoremstyle{remark}

\newcommand{\bba}{{\mathbb A}}

\newcommand{\bbc}{{\mathbb C}}

\newcommand{\bbq}{{\mathbb Q}}
\newcommand{\bbr}{{\mathbb R}}

\newcommand{\bbz}{{\mathbb Z}}

\newcommand{\Del}{{\Delta}}

\newcommand{\lam}{{\lambda}}

\newcommand{\lamb}{\underline{\lambda}}

\newcommand{\gA}{{\mathfrak A}}

\newcommand{\gB}{{\mathfrak B}}

\newcommand{\gC}{{\mathfrak C}}

\newcommand{\gM}{{\mathfrak M}}

\newcommand{\gp}{{\mathfrak p}}

\font\tenscr=rsfs10 

\newcommand{\sE}{\hbox{\tenscr E}}

\newcommand{\sS}{\hbox{\tenscr S}}

\newcommand{\cM}{{\mathcal M}}
\newcommand{\cN}{{\mathcal N}}

\newcommand{\cK}{{\mathcal K}}

\newcommand{\cT}{{\mathcal T}}
\newcommand{\cR}{{\mathcal R}}
\newcommand{\co}{{\mathcal O}}
\newcommand{\cD}{{\mathcal D}}

\newcommand{\aut}{{\operatorname{Aut}}\,}

\newcommand{\aff}{{\operatorname {Aff}}}

\newcommand{\n}{{\operatorname {N}}}

\newcommand{\gl}{{\operatorname{GL}}}

\newcommand{\spl}{{\operatorname{SL}}}

\newcommand{\sst}{{\operatorname{ss}}}

\newcommand{\id}{{\operatorname{id}}}

\newcommand{\ord}{{\operatorname{ord}}}
\newcommand{\sym}{{\operatorname{Sym}}}

\newcommand{\Hom}{{\operatorname{Hom}}}

\newcommand{\br}{{\operatorname{B}}}

\newcommand{\bs}{Schwartz--Bruhat function}
\newcommand{\psf}{Poisson summation formula}

\newcommand{\res}{\operatorname{Res}}

\newcommand{\supp}{\operatorname{supp}}

\newcommand{\A}{\bba}
\newcommand{\Z}{\bbz}
\newcommand{\Q}{\bbq}
\newcommand{\R}{\bbr}
\newcommand{\C}{\bbc}

\newcommand{\ma}{\bba^{\times}}
\newcommand{\mr}{\R_+}

\newcommand{\mk}{k^{\times}}

\newcommand{\md}{d^{\times}}
\newcommand{\gMf}{\gM_{\text f}}

\newcommand{\ti}{\widetilde}

\newcommand{\ac}[1]{\langle{#1}\rangle}

\newcommand{\twtw}[4]
{\begin{pmatrix}{#1}&{#2}\\{#3}&{#4}\\\end{pmatrix}}

\newcommand{\aaa}{\mathfrak a}
\newcommand{\ccc}{\mathfrak c}
\newcommand{\stnz}{{\operatorname{St}}}
\newcommand{\diag}{{\operatorname{diag}}}
\newcommand{\cl}{{\operatorname{Cl}}}
\newcommand{\phiu}{\Phi}
\newcommand{\phiw}{\Psi}
\newcommand{\phis}{\Upsilon}

\newcommand{\Af}{\mathbb A_{\mathrm f}}
\newcommand{\maf}{\mathbb A_{\mathrm f}^\times}
\newcommand{\z}{\overline z}
\newcommand{\y}{\overline y}
\newcommand{\gaa}{G_\aaa}
\newcommand{\vaa}{V_\aaa}
\newcommand{\vac}{V_{\aaa,\ccc}}
\newcommand{\bac}{B_{\aaa,\ccc}}
\newcommand{\wac}{W_{\aaa,\ccc}}
\newcommand{\gpX}{\mathfrak X}
\newcommand{\vcX}{\mathcal X}
\newcommand{\vcY}{\mathcal Y}
\newcommand{\vcZ}{\mathcal Z}
\newcommand{\type}{\mathcal T}
\newcommand{\pr}{\mathrm{pr}}
\newcommand{\fin}{{\rm f}}
\renewcommand{\gMf}{\gM_{\rm f}}

\begin{document}

\title[distributions of cubic algebras]
{Distributions of discriminants of cubic algebras}

\author[Takashi Taniguchi]{Takashi Taniguchi}
\address{Department of Mathematical Sciences, University of Tokyo}
\email{tani@ms.u-tokyo.ac.jp}
\date{\today}

\begin{abstract}
We study the space of binary cubic and quadratic forms
over the ring of integers $\co$ of an algebraic number field $k$.
By applying the theory of prehomogeneous vector spaces
founded by M. Sato and T. Shintani,
we can associate the zeta functions for these spaces.
Applying these zeta functions,
we derive some density theorems on
the distributions of discriminants of cubic algebras of $\co$.
In the case $k$ is a quadratic field,
we give a correction term as well as the main term.
These are generalizations of Shintani's asymptotic formulae
of the mean values of class numbers of binary cubic forms over $\mathbb Z$.
\end{abstract}

\maketitle

\section{Introduction}\label{sec:introduction}
Let $k$ be a number field and $\co$ the ring of integers of $k$.
Let $r_1$ and $r_2$ be the number of real and complex places of $k$.
We denote by $\Delta_k$, $h_k$ and $\zeta_k(s)$
the absolute discriminant, the class number
and the Dedekind zeta function of $k$, respectively.
We put
\begin{equation*}
\gA_k:=(\res_{s=1}\zeta_k(s))\cdot\frac{\zeta_k(2)}{2^{r_1+r_2+1}},\quad
\gB_k:=
	(\res_{s=1}\zeta_k(s))\cdot\frac{3^{r_1+r_2/2}\zeta_k(1/3)}{5\cdot 2^{r_1+r_2}\Delta_k^{1/2}}
\left(\frac{\Gamma(1/3)^3}{2\pi}\right)^{[k:\Q]}.
\end{equation*}

For $0\leq i\leq r_1$, let $h_i(n)$ be the numbers of the following set:
\begin{equation*}
h_i(n):=\#
\left\{(R,F)\ \vrule\begin{array}{l}
\text{$F$ is a cubic extension of $k$ with $r_1+2i$ real places,}\\
\text{$R$ is an order of $F$ containing $\co$,
and $N(\Delta_{R/\co})=n$.}\\\end{array}
\right\}.
\end{equation*}
Here $\Delta_{R/\co}$ is the relative discriminant of $R/\co$
(which is an integral ideal of $\co$)
and $N(\Delta_{R/\co})$ is its ideal norm.
Note that we count pairs $(R,F)$ up to isomorphism.
One of the primary purpose of this paper is to investigate
the function $\sum_{n<X}h_i(n)$ as $X\to\infty$.
Here we state our result when $k$ is a quadratic field.

\begin{thm}
[Theorem \ref{thm:field}]
\label{thm:introfield}
Let $k$ be a quadratic field. For any $\varepsilon>0$,
\begin{equation*}
{\tiny \begin{pmatrix}r_1\\i\\\end{pmatrix}}^{-1}
\sum_{n<X}h_i(n)
=3^{-i-r_2}\mathfrak A_kX
+3^{-i/2}\mathfrak B_kX^{5/6}
+O(X^{9/11+\varepsilon})
\qquad (X\to\infty),
\end{equation*}
where $\left(\begin{smallmatrix}r_1\\i\\\end{smallmatrix}\right)$
is the binomial coefficient.
\end{thm}

The case $k=\Q$ is known by Shintani \cite{shintanib}.
For the case $[k:\Q]\geq3$, see Theorem \ref{thm:field}.

We explain one more theorem we prove in this paper.
We call a finite $\co$-algebra a {\em cubic algebra}
if it is projective of rank $3$ as an $\co$-module.
We denote by $\mathcal C(\co)$ the set of isomorphism
classes of cubic algebras of $\co$.
For a fractional ideal $\aaa$ of $k$, we put
$\mathcal C(\co,\aaa)=\{R\in C(\co)\mid \bigwedge^3R\cong\aaa\}$.
It is known that
$\mathcal C(\co,\aaa)$ depends only on the ideal class of $\aaa$
and that $\mathcal C(\co)=\coprod_{\aaa\in\cl(k)}\mathcal C(\co,\aaa)$
(we use the same symbol $\aaa$ to denote its ideal class.)
In general for a projective $\co$-module $M$ of rank $m$,
the class of the ideal isomorphic to $\bigwedge^mM$ is called
the {\em Steinitz class} of $M$.

We count the number of $\mathcal C(\co,\aaa)$ for each $\aaa$.
More precisely, for $0\leq i\leq r_1$ we count
$\mathcal C(\co,\aaa)_i=\{R\in\mathcal C(\co,\aaa)\mid
R\otimes_\Z \R\cong \R^{r_1+2i}\times\C^{3r_2+r_1-i}\}$.
An interesting phenomenon we prove
in the case $k$ is a quadratic field is that,
the Steinitz class is
{\em not uniformly distributed} in the $X^{5/6}$-term
if $\cl(k)$ contains a non-trivial $3$-torsion element.

\begin{thm}
[Theorem \ref{thm:densityGV}]
\label{thm:introalgebra}
Let $k$ be a quadratic field. For any $\varepsilon>0$,
\vspace{-1mm}
\begin{equation*}
{\tiny \begin{pmatrix}r_1\\i\\\end{pmatrix}}^{-1}
\hspace{-5mm}
\sum_{\underset{N(\Delta_{R\slash \mathcal O})\leq X}
	{R\in \mathcal C(\mathcal O,\mathfrak a)_i}}
\frac{1}{{}^{\#}(\aut(R))}=(1+\frac1{3^{i+r_2}})\frac{\mathfrak A_k}{h_k}X
+\tau(\mathfrak a)\frac{\mathfrak B_kh_k^{(3)}}{3^{i/2}h_k}X^{5/6}
+O(X^{7/9+\varepsilon})
\quad (X\to\infty).
\vspace{-1mm}
\end{equation*}
Here ${}^{\#}(\aut(R))$ denote the cardinality
of the automorphisms of $R$ as an $\co$-algebra
and  $h_k^{(3)}$ be the number of $3$-torsions of $\cl(k)$
(which is a power of $3$.)
Also for $\mathfrak a\in \cl(k)$, we put $\tau(\mathfrak a)=1$
if there exists
$\mathfrak b\in \cl(k)$ such that $\mathfrak a=\mathfrak b^3$
and $\tau(\mathfrak a)=0$ otherwise.
\end{thm}

\medskip

Before explaining the contents of this
paper, we include a brief historical
overview of the area in order to clarify
the background of our results.
The study of distributions of discriminants or
mean values of class numbers
of rings or fields extensions is a classical topic
in algebraic number theory and has a long history.
The roots of this topic traces to Gauss,
who is the first mathematician
introducing a group theoretical approach to number theory.
In \cite{gauss},
he found that the set of orbits $\gl(2)_\Z\backslash \sym^2\Z^2$
of integral binary quadratic forms corresponds bijectively
to the set of ideal classes of {\em quadratic} rings.
Using this he gave conjectures for the asymptotic property of
the average number of class numbers of quadratic rings.
This conjecture was first proved
by Lipschitz for the imaginary case,
and by Siegel for the real case.
Siegel \cite{siegele} also proved
a density theorem for $\gl(n)_\Z\backslash \sym^2\Z^n$ in general.

If one consider {\em cubic} object as we shall do in this paper,
the most basic representation is 
the space of binary cubic forms $(\gl(2),\sym^3\aff^2)$.
The interpretation of set of orbits
of the space of integral binary cubic forms
$\gl(2)_\Z\backslash\sym^3\Z^2$ in terms of cubic rings
were discovered by Delone-Faddeev \cite{defa}
and many applications to number theory or
representation theory are obtained to the present.
For example, this was used by Davenport-Heilbronn \cite{daheb}
to prove the density of discriminants of cubic fields
$\sum_{[F:\Q]=3,|\Delta_F|<X}1\sim 3^{-1}\zeta(3)^{-1}X \ (X\to\infty)$.

In 1972, Shintani \cite{shintania}
made a significant contribution
to the study of the class numbers of $\gl(2)_\Z\backslash\sym^3\Z^2$
by applying the {\em zeta function} theory
of {\em prehomogeneous vector spaces}
founded by M. Sato and Shintani \cite{sash}.
He gave the analytic continuations, functional equations
and residue formulae of the Dirichlet series.
Combined with the results of zeta functions of binary quadratic forms
\cite{shintanib},
he gave a correction term \cite[Theorem 4]{shintanib}
to the main term of Davenport's asymptotic formula of
distributions of discriminants of irreducible cubic rings over $\Z$.

His works \cite{shintania, shintanib} are the starting point of our work.
The zeta functions of prehomogeneous vector spaces are
defined over any algebraic number fields,
and by the use of adelic language numerous contributions
to number theory has obtained.
For the case of the space of binary cubic forms,
a series of work by Wright \cite{wright} and
Datskovsky-Wright \cite{dawra, dawrb}
gave the generalization of the Davenport-Heilbronn's
density theorem over $\Q$ above
to over a general number field $k$
with finite number of local conditions.
This was recently improved by Kable-Wright \cite{kawr}
to prove an equidistribution result for
the Steinitz classes of cubic extensions.
In both cases what is called filtering process was used
to count cubic field extensions of $k$
which corresponds a set of {\em rational} equivalence classes,
rather than the cubic algebras of $\co$
which to a set of {\em integral} equivalence classes.
The density problems of integral equivalence classes
are not well considered other than over $\Z$,
and this is what we focus in this paper.
In the integral equivalence case, Landau's Tauberian theorem \cite{landaub}
modified by Sato-Shintani \cite{sash} gives a sharp error estimate
since our zeta function satisfies a functional equation.
(Note that a part of Theorem \ref{thm:densityGV}
is obtained by Wright in his thesis \cite{wrightthesis}
as the mean value of class numbers.)

We first prove Theorem \ref{thm:introalgebra} (Theorem \ref{thm:densityGV})
and after that Theorem \ref{thm:introfield} (Theorem \ref{thm:field}).
The step is to separate the reducible algebras of $\mathcal C(\co)$,
i.e., those $R\in\mathcal C(\co)$ with $R\otimes_\co k$ not a field.
As in the case of $\Z$ treated in \cite{shintanib},
these are parameterized essentially by $(\br(2),\sym^2\aff^2)$
where $\br(2)$ is the Borel subgroup of $\gl(2)$
consisting of lower triangular matrices.
We generalize the argument of \cite{shintanib} to general number fields
both algebraically and analytically.

To the author's best knowledge, Theorem \ref{thm:introalgebra}
is the first unequidstribution result of the Steinitz classes.
We mention that the idea that the character in the zeta function
segregates cubic algebras by their Steinitz classes
is due to Kable-Wright \cite{kawr}.
For a conjecture of distributions of Steinitz classes for fields,
see \cite[Introduction]{kawr}.

\medskip
After the summary of the theory of binary cubic forms,
a brief illustration of some resent progress and arising problems
of related topics is appropriate.
The parameterizations of extensions degree $2$, $3$, $4$ and $5$
was first systematically established in the celebrated
work of Wright-Yukie \cite{wryu}
via $8$ prehomogeneous representations
containing $(\gl(2),\sym^2\aff^2)$ and $(\gl(2),\sym^3\aff^2)$
as basic cases.
A program to prove a series of density theorems were proposed,
and besides, deep contributions to algebraic number theory were indicated.
After Wright-Yukie, the arithmetic theory of $\Z$-orbits was first handled
by Kable \cite{kable-3tensor} for one special case,
and recently being developed greatly by Bhargava
\cite{bha,bhb,bhc,bhe}.
In his published works,
he gave generalizations of Gauss' and Delone-Faddeev's orbit ring maps
in the degree $2$, $3$ and $4$ case,
and also obtained the main term of the density of discriminants
of quartic fields and rings.
He also announced some results for the quintic case.
On the other hand the quintic case is also handled slightly earlier
by Kable-Yukie \cite{kayua, kayub, kayuc}
and an upper bound of the number of quintic fields
ordered by the size of discriminants is obtained.

All of the $8$ prehomogeneous vector spaces treated in \cite{wryu}
are what we call as {\em parabolic type}
classified by Rubenthaler in his thesis \cite{rubea}.
These are obtained by choosing a semisimple group
and a maximal parabolic subgroup.
For details, see \cite{wryu}, \cite{yukiec} or their references also.
The case of binary cubic form corresponds to the exceptional group $G_2$
and its Heisenberg parabolic subgroup,
and our results might contribute to the theory of $G_2$ such as \cite{ggs}.
In addition, much progress of the theory of zeta function of prehomogeneous
vector space should be done both globally and locally.
One more problem we would like to mention due to its
number theoretical interest is various kinds of generalizations
of Ohno \cite{ohno} and Nakagawa's \cite{nakagawa}
extra functional equation of the global zeta function
for the space of binary cubic forms over $\Q$
to such as general number fields or other prehomogeneous representations.
We hope these theory to be developed in the future.

\medskip
We conclude this section with a brief review of contents of this paper.
In Section \ref{sec:pv} we define the representations
we consider in this paper.
Main objects we consider are the space of binary cubic forms
$(G,V)=(\gl(2),\sym^3\aff^2)$
and the space of binary quadratic forms $(B,W)$,
where $W=\sym^2\aff^2$ and $B$ is close to
the Borel subgroup of $\gl(2)$ consisting of lower triangular matrices.
An ``embedding'' of $(B,W)$ into $(G,V)$ described
in Definition \ref{defn:embedding}
plays an essential role when removing the contributions of reducible algebras.
The algebraic part of this paper is developed
in Section \ref{sec:parameterization}.
We consider the group theoretical parameterization of
cubic algebras by means of $(G,V)$ and $(B,W)$.
Since arithmetic plays no roles here,
we consider it over general Dedekind domains.
For the space of binary cubic forms $(G,V)$,
this is regarded as a generalization of
Delone-Faddeev's orbit ring map \cite{defa}.

The rest of this paper is devoted to the analytic theory.
After we introduce notation for number fields and give
normalizations of invariant measure in Section \ref{sec:measure},
we concentrate on the analysis of the zeta functions.
Both $(G,V)$ and $(B,W)$ are typical examples of
prehomogeneous vector spaces and the associated zeta functions
are studied by many mathematicians after Shintani's
pioneering works \cite{shintania,shintanib}.
As we desire to work over general number fields,
we need to rewrite his work into adelic language.
For $(G,V)$, this is done by Wright \cite{wright}.
In Section \ref{sec:global} we give the adelic version
of the zeta function of $(B,W)$ in \cite[Chapter 1]{shintanib}.
We choose F. Sato's modified approach \cite{fsatoa}
where the ``enlarged representation''
$(H,U)=(\gl(1)\times\gl(2),\sym^2\aff^2\oplus\aff^2)$ is well used.
We note that the adelic zeta function for $(H,U)$ is handled
by Yukie \cite[Chapter 7]{yukiec} in a slightly different formulation.

In Section \ref{sec:local} we deal with the archimedean local theory.
Since this determines the gamma factor
of the functional equations, this is important for our purposes.
Fortunately this is well established and we briefly recall it.
In Section \ref{sec:density} we study the density theorems.
We define our target Dirichlet series and describe them
by means of the integral expressions of the global zeta functions.
Combined with Sato-Shintani's Tauberian theorem \cite[Theorem 3]{sash},
we find the asymptotic formulae.

In Section \ref{sec:cubic}, we study some zeta integrals
for the space of binary cubic forms $(G,V)$.
Especially, the residue of the second pole
as well as the meromorphic continuation
of the zeta integral corresponds to irreducible cubic forms is obtained.
We hope this result will be useful when attacking the
conjecture of Datskovsky-Wright \cite{dawrc} or Roberts \cite{roberts}
on the second term of the asymptotics of cubic field discriminants.
We note that some of the result is obtained by Yukie
in an unpublished note \cite{yukies} by a fairly different method.

In Appendices we give some supplements to Section \ref{sec:local}
those are logically independent of our main theorems.
Especially the explicit formula of the local zeta function
for the standard test function at finit places is given.
This result is closely related to Ibukiyama-Saito's work \cite{ibsab}.

\bigskip
\noindent
{\bf Notation.}
For a finite set $X$ we denote by ${}^\# X$ its cardinality.
The standard symbols $\Q$, $\R$, $\C$ and $\Z$ will denote respectively
the set of rational, real and complex numbers and the rational integers.
If $V$ is a scheme defined over a ring $R$ and $S$ is an $R$-algebra
then $V_S$ denotes its $S$-rational points.
(We do {\em not} use the notation $V_S$ in the sense of the base change.)
If an abstract group $G$ acts on a set $X$,
then for $x\in X$ we set ${\rm Stab}(G;x)=\{g\in G\mid gx=x\}$.
If $\mathfrak x\in G\backslash X$ is the class of $x\in X$,
we also denote the group by ${\rm Stab}(G;\mathfrak x)$,
which is well defined up to isomorphism.
We always regard the affine $n$-space $\aff^n$
the set of row vectors and $\gl(n)$ acts on this space from the right.

\section{Prehomogeneous vector spaces}\label{sec:pv}
In this section, we introduce representations
$(G,V)$, $(B,W)$, $(H,W)$ and $(H,U)$
we consider in this paper
and discuss their basic properties.
The first one is the space of binary cubic forms.
The remaining three are concerning on the space of binary quadratic forms
and these are closely related.
All of these are what we call prehomogeneous vector spaces
in the sense of \cite{sash}.

\subsection{The space of binary cubic forms}
Let $V=\sym^3\aff^2$.
We regard $V$ as the space of binary cubic forms of variables
$v=(v_1,v_2)$. Elements of $V$ are expressed in the form
$x=x_0v_1^3+x_1v_1^2v_2+x_2v_1v_2^2+x_3v_2^3.$
We choose $x=(x_0,x_1,x_2,x_3)$ as the coordinate system of $V$.
We define the action of $G=\gl(2)$ on $V$ by
\begin{equation*}
(gx)(v)=(\det g)^{-1}x(vg)=\frac{1}{ad-bc}x(av_1+cv_2,bv_1+dv_2),
\quad g=\twtw abcd\in G,
\end{equation*}
which is slightly different from the usual one.
Note that this twisted representation is faithful,
whereas the usual action has kernel $\mu_3$.
Let $P(x)$ be the discriminant of $x$;
\begin{equation*}
P(x)=x_1^2x_2^2-4x_0x_2^3-4x_1^3x_3+18x_0x_1x_2x_3-27x_0^2x_3^2.
\end{equation*}
We put $\chi(g)=\det(g)$.
Then we have $P(gx)=\chi(g)^2P(x)$.
We put $V^\sst=\{x\in V\mid P(x)\not=0\}$,
which is a single $G$-orbit over any algebraically closed field.
\subsection{The spaces of binary quadratic forms and an embedding}
Let $W=\sym^2\aff^2$, regard as the space of binary quadratic forms
of variables $v=(v_1,v_2)$, and express elements
as $y=y(v)=y_1v_1^2+y_2v_1v_2+y_3v_2^2$.
We choose $y=(y_1,y_2,y_3)$ as the coordinate system of $W$.
Let $H=\gl(1)\times\gl(2)$ and express elements of $H$ as
$h=(t,g)$ where $t\in\gl(1)$ and $g\in\gl(2)$.
We define characters $\chi_1$ and $\chi_2$ on $H$
by $\chi_1(h)=t$ and $\chi_2(h)=t\det(g)$.
We define the action of $H$ on $W$ by
%
\begin{equation*}
(hy)(v)=ty(vg)=ty(av_1+cv_2,bv_1+dv_2),
\quad h=(t,g)=\left(t,\twtw abcd\right)\in H.
\end{equation*}
%
We define a subgroup $B$ of $H$ by
\begin{equation}\label{eq:B}
B=\left\{b=\left(t,\twtw 10up\right)
	\ \vrule\ t,p\in\mathbb G_m,u\in\mathbb G_a\right\}
\cong \mathbb G_m^2\ltimes \mathbb G_a,
\end{equation}
and consider the representation $(B,W)$.
We put $Q_1(y)=y_1$ and $Q_2(y)=y_2^2-4y_1y_3$.
Then we have $Q_1(by)=\chi_1(b)Q_1(y)$,
$Q_2(by)=\chi_2(b)^2Q_2(y)$ for $b\in B$, $y\in W$.
We put $W^\sst=\{y\in W\mid Q_1(y)Q_2(y)\not=0\}$,
which is a single $B$-orbit over any algebraically closed field.
We also consider the representation $(H,W)$.
This representation has a single basic relative invariant polynomial $Q_2(y)$
with the character $\chi_2^2$.
We put $W\, \tilde{}=\{y\in W\mid Q_2(y)\not=0\}$,
which is a single $H$-orbit over any algebraically closed field.

Let $S=\aff^2$. We express elements of
$S$ as $\y=(\y_1,\y_2)$, and choose this as the coordinate system of $S$.
We put $U=W\oplus S$ and express elements of $U$ as
$\ti y=(y,\y)=(y_0,y_1,y_2,\y_1,\y_2)$.
Then we can define the action of $H$ on $U$ as follows:
\begin{equation*}
h\ti y=h(y,\y)=(hy,\y g^{-1}),\qquad h=(t,g)\in H.
\end{equation*}
We put $R_1(\ti y)=y_0\y_1^2+y_1\y_1\y_2+y_2\y_2^2$ and
$R_2(\ti y)=Q_2(y)=y_2^2-4y_1y_3$.
Then we have $R_1(h\ti y)=\chi_1(h)R_1(\ti y)$
and $R_2(h\ti y)=\chi_2(h)^2R_2(\ti y)$.
We put $U^\sst=\{\ti y\in U\mid R_1(\ti y)R_2(\ti y)\not=0\}$,
which is a single $H$-orbit over any algebraically closed field.

In later sections, we use an ``embedding'' $(B,W)$
into $(G,V)$ which will play a significant role in this paper.
We define the embedding here.
\begin{defn}\label{defn:embedding}
For a binary quadratic form $y\in W$, we regard
the binary cubic form $y^\ast=v_2y$ as an element of $V$.
Then $P(y^\ast)=Q_1(y)^2Q_2(y)$.
We embed $B$ into $G$ via the map
\begin{equation*}
B\ni b=\left(t,\twtw 10us\right)\longmapsto b^\ast=\twtw t0{tu}{ts}\in G.
\end{equation*}
Then for $y\in W$ and $b\in B$, we have $(by)^\ast=b^\ast y^\ast\in V$.
\end{defn}
\begin{table}
\begin{center}
\begin{tabular}{|c|c|c|c|}
\hline
Representation & Group & Vector space & R.I. \& Characters \\
\hline\hline
$(G,V)$	&$\gl(2)$		&$\sym^3\aff^2$
	&$P\leftrightarrow \chi^2$\\
\hline\hline
$(B,W)$	&as in \eqref{eq:B}	&$\sym^2\aff^2$
	&$Q_1\leftrightarrow \chi_1,Q_2\leftrightarrow \chi_2^2$\\
\hline
$(H,W)$	&$\gl(1)\times\gl(2)$	&$\sym^2\aff^2$
	&$Q_2\leftrightarrow \chi_2^2$\\
\hline
$(H,U)=(H,W\oplus S)$
	&$\gl(1)\times\gl(2)$	&$\sym^2\aff^2\oplus\aff^2$
	&$R_1\leftrightarrow \chi_1,R_2\leftrightarrow \chi_2^2$\\
\hline
\end{tabular}
\end{center}
\vspace*{5pt}
\caption{Prehomogeneous vector spaces}%
\end{table}

\section{Parameterization of cubic algebras over a Dedekind domain}
\label{sec:parameterization}
Throughout this section, 
we assume $\mathcal O$ a Dedekind domain and $k$ its quotient field.
In this section we give group theoretical parameterizations
of $\co$-algebras which are projective of rank $3$ as $\co$-modules
by means of the representations $(G,V)$ and $(B,W)$.
The results are Propositions 
\ref{dfd} and \ref{prop:reducibleparameterization}.
The notion of Steinitz class
naturally arises in the process.
A review of geometric interpretations of orbits over fields are included.
\subsection{Projective modules over a Dedekind domain}
For the convenience of reader,
we review the basic properties of projective modules
and fractional ideals of Dedekind domains.
For details, see Milnor's book \cite[\S1]{milnorb} for example.
We assume all the $\mathcal O$-modules to be finitely generated.
\begin{defn}\label{defn:projective}
\begin{enumerate}[{\rm (1)}]
\item
An element $m$ of an $\co$-module $M$ is called a {\em torsion element}
if there exists a non-zero element $a\in\co$ such that $am=0$.
The set of torsion elements are called the {\em torsion submodule} of $M$.
If the torsion submodule of $M$ is trivial, $M$ is called {\em torsion free}.
\item
For an $\co$-module $M$, $M\otimes k$ is a vector space over $k$.
The dimension is called the {\em rank} of $M$.
\end{enumerate}
\end{defn}
\begin{prop}\label{prop:projective}
\begin{enumerate}[{\rm (1)}]
\item
If $M$ is torsion free, then the map $M\rightarrow M\otimes k$ is injective.
\item
An $\co$-module $M$ is torsion free if and only if $M$ is projective.
Especially, any fractional ideal of $\co$ is projective,
and any submodule of a projective module is projective.
\item
Let $M$ be a projective $\co$-module of rank $n$.
Then there exists a non-zero ideal $\aaa$ of $\co$ such that
$M\cong \co^{n-1}\oplus\aaa$ as $\co$-modules.
The ideal class of $\aaa$ is uniquely determined by $M$ and
called the {\em Steinitz class} of $M$, which we denote by $\stnz(M)$.
\item
For projective modules $M_1$ and $M_2$,
we have $\stnz(M_1\oplus M_2)=\stnz(M_1\otimes M_2)=\stnz(M_1)\cdot\stnz(M_2)$
where the last product is of the ideal class group.
\item
Let $\aaa$ be a fractional ideal of $\co$.
Then the inverse ideal $\aaa^{-1}$ of $\aaa$ is given by
$\{x\in k\mid x\aaa\subset\co\}$.
Moreover, $\Hom_{\co\text{-module}}(\aaa,\co)\cong \aaa^{-1}$
as $\co$-modules.
\item
For any non-zero fractional ideals $\aaa,\mathfrak b$ and $\ccc$,
there exist elements $x,y\in k$ such that $x\aaa+y\mathfrak b=\ccc$.
\end{enumerate}
\end{prop}

\subsection{A generalization of Delone-Faddeev's orbit ring map}
Let $\mathfrak a$ be a fractional ideal of $\mathcal O$.
We use the same symbol $\mathfrak a$ to denote its ideal class
if there is no confusion. We denote by $\cl(k)$ the ideal class group of $k$.
\begin{defn}
Let $\mathcal C(\mathcal O)$ be the set of isomorphism classes
of finite $\mathcal O$-algebras which are projective of rank 3
as $\mathcal O$-modules.
Elements of $\mathcal C(\mathcal O)$ are called
{\em cubic algebras}. For any fractional ideal $\mathfrak a$,
we define
$\mathcal C(\mathcal O,\mathfrak a)=
	\{R\in\mathcal C(\mathcal O)\mid\stnz(R)=\aaa\}$.
\end{defn}
\begin{prop}\label{prop:unit}
\begin{enumerate}[{\rm (1)}]
\item
We have
$\mathcal C(\mathcal O)
	=\coprod_{\mathfrak a\in \cl(k)}\mathcal C(\mathcal O,\mathfrak a)$.
\item
Let $R\in\mathcal C(\mathcal O,\mathfrak a)$.
There is an isomorphism $R\slash\mathcal O\cong\mathcal O\oplus \mathfrak a$
of $\mathcal O$-modules.
\end{enumerate}
\end{prop}
\begin{proof}
(1) immediately follows from Proposition \ref{prop:projective} and
we consider (2).
We prove $R/\co$ is torsion free.
Let $r\in R$ and non-zero $x\in\co$ satisfy $xr\in\co$.
Then $r\in k$. Since $R$ is a finite $\co$-algebra,
$r$ is integral over $\co$. Hence we have $r\in\co$
because $\co$ is integrally closed in $k$.
This shows that $R/\co$ is torsion free and hence projective.
We have $R\cong \co\oplus(R/\co)$ and
therefore $R/\co\cong \co\oplus\aaa$ by Proposition \ref{prop:projective} (3).
\end{proof}
We consider a parameterization of 
$\mathcal C(\mathcal O,\mathfrak a)$ using the representation $(G,V)$,
which is a generalization of Delone-Faddeev's orbit ring map
\cite[Section 15]{defa} over $\Z$.

\begin{defn}We put
\begin{align*}
G_k\supset \gaa&=
\left\{\twtw abcd\ \vrule \ a\in\co,b\in\aaa,c\in\aaa^{-1},
d\in\co,ad-bc\in\co^\times\right\},\\
V_k\supset \vaa&=
\{x\mid x_0\in \aaa, x_1\in \mathcal O,
x_2\in\aaa^{-1},x_3\in\aaa^{-2}\}.
\end{align*}
Then $\gaa\cdot \vaa\subset \vaa$.
\end{defn}
We could naturally construct isomorphisms
$\gaa\cong \aut(\co\oplus\aaa)$ and
$\vaa\cong\sym^3(\co\oplus\aaa)\otimes\wedge^2(\co\oplus\aaa)$,
such that the canonical action compatible.
\begin{prop}\label{dfd}
\begin{enumerate}[{\rm (1)}]
\item
There exists the canonical bijection
between ${\mathcal C}(\co,\aaa)$ and $\gaa\backslash \vaa$
making the following diagram commutative:
\begin{equation*}
\xymatrix{
\gaa\backslash \vaa
\ar[r]^{\ \ \ \quad\quad\ \ \ } \ar[d]^{P}
& {\mathcal C}(\co,\aaa) 
\ar[d]^{\text{discriminant}}\\
(\co^\times)^2\backslash \aaa^{-2}
\ar[r]^{\times\aaa^2\quad\quad}
&
\{\text{integral ideals of $\co$}\}.
}
\end{equation*}
Here, the right vertical arrow is to take the discriminant,
and the low horizontal arrow is given by multiplying $\aaa^2$.
Moreover, this diagram is functorial with respect to
ring homomorphisms of Dedekind domains.
\item
For each $R\in\mathcal C(\co,\aaa)$,
let $x_R$ be the corresponding element in $\gaa\backslash \vaa$.
Then $\aut_{\co\textrm{-alg}}(R)\cong \mathrm{Stab}(\gaa;x_R)$.
\end{enumerate}
\end{prop}
\begin{proof}
Since the proof is similar to \cite[Section 15]{defa}
or \cite[Proposition 4.2]{ggs}, we shall be brief.
For each $R\in\mathcal C(\co,\aaa)$, the binary cubic form
\begin{equation*}
x_R\colon
R/\co\longrightarrow \wedge^2(R/\co),
\qquad
\xi\longmapsto \xi\wedge\xi^{2}
\end{equation*}
can be regarded as an element of $\gaa\backslash \vaa$
since $R/\co\cong \co\oplus\aaa$.
This map $R\mapsto x_R$ gives the desired bijection.
To see this map in fact bijective,
we will write down this correspondence explicitly.

Let $R\in\mathcal C(\co,\aaa)$.
By Proposition \ref{prop:unit}
we fix an $\co$-module isomorphism $R\cong \co\oplus\co\oplus\aaa$
such that $(1,0,0)$ is the multiplicative identity $1$ of $R$.
We regard $R$ as a subalgebra of $R\otimes k\cong k\oplus k\oplus k$.
Let $\omega_1=(0,1,0), \omega_2=(0,0,1)$.
Then $R=\{p+q\omega_1+r\omega_2\mid p,q\in\co, r\in\aaa\}$.
Let $\omega_1\omega_2=\alpha+\beta\omega_1+\gamma\omega_2$.
Then since $\omega_1(\aaa\omega_2)\subset R$, we have
$\alpha,\beta\in\aaa^{-1},\gamma\in\co$.
Hence by replacing $\omega_1, \omega_2$ to
$\omega_1-\gamma,\omega_2-\beta$ if necessary,
we may assume $\omega_1\omega_2\in\aaa^{-1}$.
Let
\begin{equation}\label{eq:composition}
\omega_1^2=j-b\omega_1+a\omega_2,\quad
\omega_2^2=l-d\omega_1+c\omega_2,\quad
\omega_1\omega_2=m.
\end{equation}
Then the associativity of the product of the algebra $R\otimes k$ requires
\begin{equation}\label{eq:klm}
j=-ac,\quad l=-bd,\quad m=-ad.
\end{equation}
Also since $\co\omega_1^2,\aaa^2\omega_2^2\subset R$, we have
\begin{equation}\label{eq:abcd}
a\in\aaa,\quad
b\in\co,\quad
c\in\aaa^{-1},\quad
d\in\aaa^{-2}.
\end{equation}
On the other side, For any $a,b,c,d$ satisfying \eqref{eq:abcd},
the $\co$-module $\co\oplus\co\omega_1\oplus\aaa\omega_2$
becomes an $\co$-algebra if we take $(1,0,0)$ as its multiplicative identity
and define the multiplication law by
\eqref{eq:composition} with \eqref{eq:klm}.
Let $\xi(v)=v_1\omega_1+v_2\omega_2$. Then by computation we have
\begin{equation}\label{eq:xiwedge}
1\wedge\xi(v)\wedge\xi(v)^2=(av_1^3+bv_1^2v_2+cv_1v_2^2+dv_2^3)\cdot
1\wedge\omega_1\wedge\omega_2.
\end{equation}
This shows that $x_R$ is a class of
$av_1^3+bv_1^2v_2+cv_1v_2^2+dv_2^3\in \vaa$.
These show that the upper horizontal arrow bijective.

To see the commutativity of the diagram,
it is enough to check locally, i.e.,
we may assume $\co$ a discrete valuation ring.
Since this compatibility is known for PID in \cite{defa}
(or by a simple computation), we have (1).

We consider (2).
Any $\ti\psi\in\aut_{\co\text{-alg}}(R)$
obviously induce  $\psi\in\aut_{\co\text{-module}}(R/\co)$.
We fix an isomorphism $R/\co\cong \co\oplus\aaa$
and regard $\psi\in\gaa$, $x_R(v)\in\vaa$.
We also define $\omega_1,\omega_2,a,b,c,d$ as above
by using this isomorphism, so that
$x_R(v)=av_1^3+bv_1^2v_2+cv_1v_2^2+dv_2^3\in \vaa$.
Let
\begin{equation*}
\omega_1'=\ti\psi(\omega_1),\quad
\omega_2'=\ti\psi(\omega_2),\quad
\text{and}\quad
\xi'(v)=v_1\omega_1'+v_2\omega_2'.
\end{equation*}
Then since $\ti\psi\in\aut_{\co\text{-alg}}(R)$,
equation \eqref{eq:composition} also holds for
the pair $(\omega_1',\omega_2')$ and hence we have
\begin{equation}\label{eq:xi'wedge}
1\wedge\xi'(v)\wedge\xi'(v)^2=(av_1^3+bv_1^2v_2+cv_1v_2^2+dv_2^3)\cdot
1\wedge\omega_1'\wedge\omega_2'.
\end{equation}
By \eqref{eq:xiwedge}, \eqref{eq:xi'wedge} and
$1\wedge\omega_1'\wedge\omega_2'
	=(\det \psi)\cdot1\wedge\omega_1\wedge\omega_2$,
we see $\psi(x_R)=x_R$ by
\begin{multline*}
\psi(x_R(v))\cdot1\wedge\omega_1'\wedge\omega_2'
=	x_R(v\psi)\cdot1\wedge\omega_1\wedge\omega_2
=	1\wedge\xi(v\psi)\wedge\xi(v\psi)^2\\
=	1\wedge\xi'(v)\wedge\xi'(v)^2
=	x_R(v)\cdot1\wedge\omega_1'\wedge\omega_2'.
\end{multline*}
Hence we have a map
$\aut_{\co\text{-alg}}(R)\to \mathrm{Stab}(\gaa;x_R)$.
To see this map an isomorphism is a routine task
and we omit the details here.
\end{proof}
\subsection{Parameterization of reducible algebras}
We next consider the parameterization of the {\em reducible algebras},
that is, those $R\in\mathcal C(\co,\aaa)$ with $R\otimes k$ not fields.
\begin{defn}
We define subsets of $\mathcal C(\co,\aaa)$ as follows:
\begin{align*}
\mathcal C(\co,\aaa)^3
	&=\{R\mid R\otimes k \text{ is a cubic field extension of $k$}\},\\
\mathcal C(\co,\aaa)^2
	&=\{R\mid R\otimes k\cong k\times F
		\text{ where $F$ is a quadratic field extension of $k$}\},\\
\mathcal C(\co,\aaa)^1
	&=\{R\mid R\otimes k\cong k\times k\times k\},\\
\mathcal C(\co,\aaa)^0
	&=\{R\mid R\otimes k \text{ is an inseparable algebra}\}.
\end{align*}
We also define subsets of $\vaa$ as follows:
\begin{align*}
\vaa^3&=\{x\mid \text{$x(v)$ is irreducible over $k$}\},\\
\vaa^2&=\{x\mid \text{$x(v)$ decomposes into degree $1$ and $2$
	irreducible polynomials over $k$}\},\\
\vaa^1&=\{x\mid \text{$x(v)$ has three distinct rational roots
	in $\mathbb P^1$}\},\\
\vaa^0&=\{x\mid \text{$x(v)$ has a multiple root in $\mathbb P^1$}\}.
\end{align*}
Finally, we put
$\mathcal C(\co,\aaa)^{\rm red}=\coprod_{0\leq i\leq 2}\mathcal C(\co,\aaa)^i$,
$\vaa^{\rm red}=\coprod_{0\leq i\leq 2} \vaa^i$.
\end{defn}
We omit the routine proof of the following lemma.
\begin{lem}
\begin{enumerate}[{\rm (1)}]
\item We have
$\mathcal C(\co,\aaa)=\coprod_{0\leq i\leq 3}\mathcal C(\co,\aaa)^i$,
$\vaa=\coprod_{0\leq i\leq 3} \vaa^i$.
\item Each $\vaa^i$ is a $\gaa$-invariant subset.
\item The set $\mathcal C(\co,\aaa)^i$ corresponds bijectively
to the set $\gaa\backslash \vaa^i$
via the map of Proposition \ref{dfd}.
\end{enumerate}
\end{lem}
Hence, we consider $\gaa\backslash \vaa^{\rm red}$.
Let $\aaa,\ccc$ be fractional ideals of $\co$.
\begin{defn}We put
\begin{align*}
B_k\supset \bac&=
\left\{\left(t,\twtw 10up\right)\ \vrule\ 
	t,p\in\co^\times,u\in\aaa^{-1}\ccc^{-2}\right\},\\
W_k\supset \wac&=
\{y\mid y_1\in\ccc, y_2\in\aaa^{-1}\ccc^{-1},y_3\in\aaa^{-2}\ccc^{-3}\}.
\end{align*}
Then $\bac\cdot \wac\subset \wac$.
\end{defn}
\begin{defn}
We define subsets of $\wac$ as follows:
\begin{align*}
\wac^2&=\{y\mid \text{$y(v)$ is irreducible over $k$}\},\\
\wac^1&=\{y\mid \text{$y(v)$ has two distinct rational roots
	in $\mathbb P^1$}\},\\
\wac^0
&=\{y\mid \text{$y(v)$ has a multiple root in $\mathbb P^1$}\}.
\end{align*}
Then $\wac=\coprod_{0\leq i\leq 2}\wac^i$
and each $\wac^i$ is a $\bac$-invariant subset.
\end{defn}
To construct the orbit ring map, we need a lemma.
\begin{lem}\label{lem:qsmn}
Let $q,s,m,n\in k$ satisfy $q\aaa^{-1}+s\co=m\aaa^{-1}+n\co=\ccc$.
Then there exists $\gamma\in\gaa$ such that
\begin{equation*}
\gamma\begin{pmatrix}m\\ n\end{pmatrix}
=\begin{pmatrix}q\\ s\end{pmatrix}.
\end{equation*}
\end{lem}
\begin{proof}
Let
$\pi\colon \aaa^{-1}\oplus\co\ni x=(x_1,x_2)\mapsto qx_1+sx_2\in\ccc$.
This is surjective, and since $\ccc$ is projective
there is a section $\mu$ of $\pi$.
Hence we have an isomorphism
\begin{equation*}
\phi\colon \aaa^{-1}\oplus\co\longrightarrow \ccc\oplus \ker\pi,
\quad
x\longmapsto (\pi(x),x-\mu\circ\pi(x)).
\end{equation*}
For
$\pi'\colon \aaa^{-1}\oplus\co\ni (x_1,x_2)\mapsto mx_1+nx_2\in\ccc$
we also have an isomorphism
\begin{equation*}
\phi'\colon \aaa^{-1}\oplus\co\longrightarrow \ccc\oplus \ker\pi',
\quad
x\longmapsto (\pi'(x),x-\mu'\circ\pi'(x)).
\end{equation*}
By Proposition \ref{prop:projective}, we have an $\co$-module isomorphism
$\varphi\colon \ker\pi\rightarrow \ker\pi'$.
Let $\Phi=\phi'^{-1}\circ(\id\oplus\varphi)\circ\phi$,
which is an element of
$\aut_{\co\text{-module}}(\aaa^{-1}\oplus\co)$;
\begin{equation*}
\xymatrix{
\aaa^{-1}\oplus\co
	\ar[r]^{\phi} \ar[d]_{\Phi} &
\ccc\oplus\ker\pi\ar[d]^{\id\oplus\varphi}\\
\aaa^{-1}\oplus\co
	\ar[r]^{\phi'}&
\ccc\oplus\ker\pi'
}
\end{equation*}
Then $\Phi$ is represented by an element $\gamma\in\gaa$
i.e., $\Phi(x)=x\gamma$ for $x\in \aaa^{-1}\oplus\co$.
By the commutativity of the diagram above,
we have $\pi'\circ\Phi(x)=\pi(x)$, i.e.,
$x\gamma \left(\begin{smallmatrix}m\\n\\\end{smallmatrix}\right)
=x\left(\begin{smallmatrix}q\\s\\\end{smallmatrix}\right)$.
Since this holds for any $x$, 
we have the lemma.
\end{proof}
\begin{prop}\label{prop:reducibleparameterization}
\begin{enumerate}[{\rm (1)}]
\item
We fix $\aaa$.
For each $\ccc$, there exists the canonical map
$\psi_{\aaa,\ccc}\colon
	\bac\backslash \wac\rightarrow
	\gaa\backslash \vaa^{\rm red}$
making the following diagram commutative:
\begin{equation*}
\xymatrix{
\bac\backslash \wac
\ar[r]^{\psi_{\aaa,\ccc}} \ar[d]_{Q_1^2Q_2} &
	\gaa\backslash \vaa^{\rm red}\ar[d]^{P}\\
(\co^\times)^2\backslash \aaa^{-2} \ar@{=}[r] &
(\co^\times)^2\backslash \aaa^{-2}.
}
\end{equation*}
\item
Moreover, collecting $\psi_{\aaa,\ccc}$ for all $\ccc\in \cl(k)$
gives the map
\begin{equation*}
\coprod_{\ccc\in \cl(k)}
	\bac\backslash \wac^2\longrightarrow
	\gaa\backslash \vaa^2
\end{equation*}
bijective.
Also for $y\in \wac^2$ we have
$\mathrm{Stab}(\bac;y)
	\cong\mathrm{Stab}(\gaa;\psi_{\aaa,\ccc}(y))$.
\end{enumerate}
\end{prop}
\begin{proof}
We fix $\aaa$. Let
\begin{align*}
\vac^{\rm red}
&=	\{(mv_1+nv_2)(lv_1^2+l'v_1v_2+l''v_2^2)\in
		\vaa^{\rm red}\mid
	m,n,l,l',l''\in k,
	m\aaa^{-1}+n\co=\ccc\},\\
\vac^2
&=	\vac^{\rm red}\cap \vaa^2,
\end{align*}
which depend only on the ideal class of $\ccc$.
Then we have
$\vaa^{\rm red}=\bigcup_{\ccc\in \cl(k)}\vac^{\rm red}$
and 
$\vaa^2=\coprod_{\ccc\in \cl(k)}\vac^2$.
(Note that the second union is disjoint
while the first one is in general not.)
We use the embedding of Definition \ref{defn:embedding}
to construct $\psi_{\aaa,\ccc}$.
We fix $q,s\in k$ such that $q\aaa^{-1}+s\co=\ccc$.
Then, $q\in\aaa\ccc$, $s\in\ccc$, and also
there exist $p\in\ccc^{-1},r\in\aaa^{-1}\ccc^{-1}$ such that
$ps-qr\in\co^\times$.
We put $g_{\aaa,\ccc}=
\left(\begin{smallmatrix} p&q\\ r&s\\\end{smallmatrix}\right)\in G_k$.
We define
\begin{equation*}
\ti\psi_{\aaa,\ccc}\colon
	\wac\longrightarrow
	\vac^{\rm red},
\qquad	y\longmapsto g_{\aaa,\ccc}y^\ast.
\end{equation*}
Then $\ti\psi_{\aaa,\ccc}(\wac^2)\subset \vac^2$
and
\begin{equation*}
P\left(\ti\psi_{\aaa,\ccc}(y)\right)
=	P(g_{\aaa,\ccc}y^\ast)
=	\det(g_{\aaa,\ccc})^2P(y^\ast)
=	(ps-qr)^2Q_1(y)^2Q_2(y).
\end{equation*}
Also considering $B$ as a subgroup of $G$ via the embedding
of Definition \ref{defn:embedding},
we see
\begin{equation}\label{eqref:intsects}
g_{\aaa,\ccc}^{-1}\gaa g_{\aaa,\ccc}\cap B_k=\bac.
\end{equation}
This shows that the map $\ti\psi_{\aaa,\ccc}$ induces
a well defined map
$\psi_{\aaa,\ccc}\colon
\bac\backslash \wac
\to \gaa\backslash \vac^{\rm red}$
making the diagram in the proposition commutative.
Let
$g_{\aaa,\ccc}'=
\left(\begin{smallmatrix} p'&q'\\ r'&s'\\\end{smallmatrix}\right)\in G_k$
where $q'\aaa^{-1}+s'\co=\ccc$,
$p'\in\ccc^{-1},r'\in\aaa^{-1}\ccc^{-1}$ and
$p's'-q'r'\in\co^\times$.
Then since
\begin{equation*}
g_{a,c}'g_{a,c}^{-1}
=\frac{1}{ps-qr}\twtw {p's-q'r}{q'p-p'q}{sr'-rs'}{s'p-r'q}\in\gaa,
\end{equation*}
we have
$g_{\aaa,\ccc}'y^\ast\in\gaa g_{\aaa,\ccc}y^\ast$.
This shows that $\psi_{\aaa,\ccc}$ does not depend on the choice of
$g_{\aaa,\ccc}$.

We claim that $\ti\psi_{\aaa,\ccc}$ is surjective. To see this, let
\begin{equation*}
x=(mv_1+nv_2)(lv_1^2+l'v_1v_2+l''v_2^2)\in\vac^{\rm red},
\quad m\aaa^{-1}+n\co=\ccc.
\end{equation*}
We take $\gamma\in\gaa$ as in Lemma \ref{lem:qsmn}.
By changing $x$ to $\gamma x$ if necessary,
we assume $m=q, n=s$. Also by a variation of Gauss' lemma,
we have $l\in\ccc^{-1},l'\in\aaa^{-1}\ccc^{-1}$, and
$l''\in \aaa^{-2}\ccc^{-1}$. 
This shows $\ti\psi_{\aaa,\ccc}$ is surjective and hence
$\psi_{\aaa,\ccc}\colon
\bac\backslash \wac
\to \gaa\backslash \vac^{\rm red}$
also.
These shows (1) and the surjectivity of the map of (2).

Let $y,y'\in \wac^2$.
We assume $\psi_{\aaa,\ccc}(y)$ and $\psi_{\aaa,\ccc}(y')$
lie in the same $\gaa$-orbit.
Then there exists $\gamma\in\gaa$ such that
$v_2y'=(g_{\aaa,\ccc}^{-1}\gamma g_{\aaa,\ccc})(v_2y)\in V_k$.
Since $y$ and $y'$ are irreducible quadratic forms,
$g_{\aaa,\ccc}^{-1}\gamma g_{\aaa,\ccc}$ must
fix the linear form $v_2$,
i.e., $g_{\aaa,\ccc}^{-1}\gamma g_{\aaa,\ccc}\in B_k$.
Combined with \eqref{eqref:intsects} we have
$b=g_{\aaa,\ccc}^{-1}\gamma g_{\aaa,\ccc}\in\bac$.
Hence $y'=by$ and therefore $y$ and $y'$ lie in the same $\bac$-orbit.
This shows that the map of (2) injective.
The second statement of (2) can be proved similarly.
\end{proof}

\begin{rem}
In \cite{shintanib}, the corresponding statement for $\Z$
is shown in the proof of Lemma 12.
\end{rem}
\subsection{Geometric interpretation of orbits over a field}
The sets of non-singular orbits of representations we defined
in Section \ref{sec:pv} have
well known geometric interpretations over a field.
We recall those facts here.
In this subsection let $\co=k$ be an arbitrary field.
\begin{defn}
We denote by $\sE_2(k)$ and $\sE_3(k)$ the sets of 
isomorphism classes of separable quadratic and cubic algebras of $k$,
respectively.
\end{defn}
We first consider the space of binary cubic forms $(G,V)$.
\begin{defn}\label{def:rod3}
For $x=x(v_1,v_2)\in V_k^\sst$, we define
\begin{align*}
Z_x	&=\mathrm{Proj}\, k[v_1,v_2]/(x(v_1,v_2)),\\
k(x)	&=\Gamma(Z_{x},\co_{Z_{x}}).
\end{align*}
We regard $k(x)$ as an element of $\sE_3(k)$,
which is possible because
elements of $V_k^\sst$ are separable cubic polynomials.
\end{defn}

\begin{prop}\label{prop:rod3}
The map $x\mapsto k(x)$
gives a bijection between
$G_k\backslash V_{k}^\sst$
and $\sE_3(k)$.
Moreover this map coincides with the
Delone-Faddeev map we defined in Proposition \ref{dfd}
in the domain of the definition.
\end{prop}

We next consider the spaces of binary quadratic forms
$(B,W)$, $(H,W)$ and $(H,U)$.
\begin{defn}
For $\ti y=(y,\y)\in U_k^\sst$, we define
\begin{align*}
Z_{\ti y}=Z_{y}	&=\mathrm{Proj}\, k[v_1,v_2]/(y(v_1,v_2)),\\
\ti k(\ti y)=\ti k(y)	&=\Gamma(Z_{y},\co_{Z_{y}}).
\end{align*}
We regard $\ti k(\ti y)=\ti k(y)$ as an element of $\sE_2(k)$.
\end{defn}\label{def:rod2}

\begin{prop}\label{prop:rod2}
\begin{enumerate}[{\rm (1)}]
\item
The map $y\mapsto \ti k(y)$
gives bijections between
$B_k\backslash W_{k}^\sst, H_k\backslash W_{k}\tilde{},$
and $\sE_2(k)$.
Also the map $\ti y\mapsto \ti k(\ti y)$
gives a bijection between
$H_k\backslash U_k^\sst$ and $\sE_2(k)$.
\item
For $y\in W_k^\sst$, we have $k(y^\ast)\cong \ti k(y)\times k$
where $y^\ast\in V_k^\sst$ is defined in Definition \ref{defn:embedding}.
\end{enumerate}
\end{prop}

\section{Notation for number fields and invariant measures}\label{sec:measure}

For the rest of this paper, we assume $k$ is a number field
and $\co$ the ring of integers of $k$.
In this section we prepare notation for number fields, and
fix various invariant measure both locally and globally.
Let $\mr=\{t\in\R^\times\mid t>0\}$.
For an integral ideal $I$ of $\co$
the ideal norm of $I$ is denoted by $N(I)$.
We extend it to general fractional ideals in the obvious manner.
Let $\gM$,
$\gM_{\infty}$, $\gM_{\text{f}}$, 
$\gM_{\R}$ and $\gM_{\C}$
denote respectively the set of all places of $k$, all infinite
places, all finite places, all real places and all complex places.
For $v\in\gM$, $k_v$
denotes the completion of $k$ at $v$ and $|\;|_v$ 
denotes the normalized
absolute value on $k_v$. If $v\in\gMf$ then $\co_v$ denotes
the ring of integers of $k_v$,
$\gp_v$ the maximal ideal of $\co_v$ and $q_v$ the cardinality of
$\co_v/\gp_v$. 
For $t\in\ k_v^\times$, we define $\ord_v(t)$ such that
$|t|_v=q_v^{-\ord_v(t)}$.
For any separable quadratic algebra $L_v$ of $k_v$,
let $\co_{L_v}$ denote the ring of integral elements of $L_v$.
That is, if $L_v$ is a quadratic extension then $\co_{L_v}$ is
the integer ring of $L_v$ and 
if $L_v=k_v\times k_v$ then $\co_{L_v}=\co_v\times\co_v$.
If $k_1/k_2$ is a finite extension of either number fields or
non-Archimedean local fields then we shall write
$\Del_{k_1/k_2}$ for the relative discriminant of the extension;
it is an ideal in the ring of integers of $k_2$.
For conventions, we let $\Del_{k_2\times k_2/k_2}$ be
the integer ring of $k_2$.
To ease the notational burden we shall use the
same symbol, $\Del_{k_1}$, for 
the classical absolute discriminant of $k_1$ over $\Q$.
Since this number generates the ideal
$\Del_{k_1}$, the resulting notational identification is harmless.
We put $\Gamma_\R(s)=\pi^{-s/2}\Gamma(s/2)$
and $\Gamma_\C(s)=(2\pi)^{1-s}\Gamma(s)$.

Returning to $k$, we let $r_1$, $r_2$, $h_k$, $R_k$ and $e_k$ be
respectively the number of real places, the number of complex
places, the class number, the regulator and the number of roots of
unity of $k$. It will be convenient to set
$\gC_k=2^{r_1}(2\pi)^{r_2}h_kR_ke_k^{-1}$.
Let $\zeta_k(s)$ be the Dedekind zeta function of $k$.

We refer to \cite{weilc} as the basic reference for fundamental properties
on adeles.
The rings of adeles and finite adeles
are denoted by $\A$ and $\Af$,
and the groups of ideles and finite ideles of $k$
are by $\ma$ and $\maf$, respectively.
We put $k_\infty=k\otimes_\Q\R$ and $\widehat\co=\co\otimes_\Z\widehat\Z$.
Note that $\widehat\co=\prod_{v\in\gMf}\co_v$,
$\Af=\widehat\co\otimes_\co k$,
and $\A=k_\infty\times \Af$.
Following as usual,
the $k_\infty$-rational points $X_{k_\infty}$
of a variety $X$ of $k$ is abbreviated to $X_\infty$.
The adelic absolute value $|\;|$ on $\ma$ is normalized so that,
for $t\in\ma$, $|t|$ is the module of multiplication by $t$
with respect to any Haar measure $dx$ on $\A$, i.e. $|t|=d(tx)/dx$.
We define $|\;|_\infty$ on $k_\infty^\times$ similarly.
Let $\A^1=\{t\in\ma\mid |t|=1\}$.
Let $\A^0$ be the unique maximal compact subgroup of $\ma$.
Suppose $[k:\Q]=n$.
For $\lam\in\R_+$, $\lamb\in \ma$ 
is the idele whose component at any infinite place is $\lam^{1/n}$
and whose component at any finite place is $1$.  
Then we have $|\lamb|=\lam$.
Let $\mathfrak d\in\ma_\fin$ be a differential idele of $k$.

From now on we give normalizations of measures.
We first prepare common notation for products of local measures.
Let $X$ be an algebraic group over $k$.
Once we normalize a local measure
$dx_v$ on $X_{k_v}$ for each $v\in\gM$,
then we always let
$dx_\infty	=\prod_{v\in\gM_\infty}dx_v$,
$dx_\fin	=\prod_{v\in\gMf}dx_v$ and
$d_\pr x	=\prod_{v\in\gM}dx_v$,
which are measures on $X_\infty$, $X_{\Af}$, and $X_\A$, respectively.
Hence $d_\pr x=dx_\infty dx_\fin$.
A global measure usually defined in a different way is denoted by 
such as $dx$. We give the ratio of $dx$ and $d_\pr x$.

For any  $v\in \gM_{\text{f}}$, 
we choose a Haar measure $dx_v$ on $k_v$ 
to satisfy $\int_{\co_v}dx_v=1$.  
We write $dx_v$ for the ordinary Lebesgue measure if $v$ is real, 
and for twice the Lebesgue measure if $v$ is imaginary.  
On the other side,
we choose a Haar measure  $dx$ on $\A$ such that $\int_{\A/k}dx=1$.
Then $d_\pr x =\Delta_k^{1/2}dx$
(see \cite{weilc}, p. 91).
For a vector space $V$, we always choose the Haar measure on $V_\A$
such that the volume of $V_\A/V_k$ is one.
If an isomorphism $V\cong \aff^n$ is fixed,
using this identification we choose the local measure
$dx_v$ on $V_{k_v}$ as the $n$-product of the measure on $k_v$
normalized above. Hence $d_\pr x=\Delta_k^{n/2}dx$.
For any $v\in \gM_{\text{f}}$, we normalize the Haar measure $\md t_v$
on $\mk_v$ such that $\int_{\co_v^{\times}} \md t_v = 1$.  
Let $\md t_v(x)=|x|_v^{-1}d x_v$ if $v\in\gM_\infty$.
We choose a Haar measure $\md t^1$ on $\A^1$ such that
$\int_{\A^1/\mk}\md t^1=1$.
Using this measure and the decomposition
\begin{equation*}
\mr\times\A^1\cong \ma,\qquad
(\lam, t^1)\mapsto \lamb t^1,
\end{equation*}
we define an invariant measure $\md t$ on $\ma$
by $\md t=\md\lam \md t^1$ where $\md\lam=d\lam/\lam$.
Then $\md_\pr t=\gC_k\md t$
(see \cite{weilc}, p. 95).
We also choose an invariant measure $\md t^0$ on $\A^0$
such that $\int_{\A^0}\md t^0=1$.

Let $\Omega^1$ be the group of characters on $\A^1/\mk$.
For $\omega\in\Omega^1$, we put $\delta(\omega)=1$ if $\omega$ is trivial
and $0$ otherwise.
We extend $\omega$ to a character on $\ma/\mk$ by assuming that
$\omega(\lamb)=1$ for any $\lamb\in\R_+$.
Let $\Omega_v$ be the group of characters on $k_v$.
If $v\in\gMf$, we put $\delta(\omega_v)=1$
if $\omega_v$ is trivial on $\co_v^\times$ and $0$ otherwise.
For a vector space $V$,
Let $\sS(V_{\A})$, $\sS(V_{k_v})$, $\sS(V_\infty)$ and $\sS(V_{\Af})$
be the spaces of \bs s on each of the indicated domains.

Let $\mathrm T(2)\subset \gl(2)$ be the set of diagonal matrices
and $\mathrm N(2)\subset \gl(2)$ the set of the
lower triangular matrices with diagonal entries $1$.
Then, $\mathrm{B(2)=T(2)\ltimes N(2)}$ is a Borel subgroup of $\gl(2)$.
(This should not be confused to the group
$B$ we defined in Section \ref{sec:pv}.)
We express
\begin{equation*}
\diag(t_1,t_2)=\twtw{t_1}00{t_2}\in\mathrm T(2),\quad
n(u)=\twtw 10u1\in\mathrm N(2).
\end{equation*}
We briefly review the standard Haar measure on 
global and local measure on $\gl(2)$.
Let $\cK(2)=\prod_{v\in\gM}\cK_v(2)$ where
$\cK_v(2)={\rm O}(2),{\rm U}(2),\gl(2)_{\co_v}$ for
$v\in\gM_\R,\gM_\C,\gMf$, respectively.
We choose an invariant measure $d\kappa, d\kappa_v$ on $\cK(2), \cK_v(2)$
such that $\int_{\cK(2)}d\kappa=1, \int_{\cK_v(2)}d\kappa_v=1$, respectively.
Obviously, $d\kappa=d_\pr\kappa$.
We express element $b\in{\rm B}(2)_\A$ as
$b=tn(u)=\diag(t_1,t_2)n(u)$ where $t_1,t_2\in\ma$ and $u\in\A$.
We put $\md t=\md t_1\md t_2$ which is a Haar measure on ${\rm T}(2)_\A$,
and choose $db=|t_2/t_1|\md t du$
as a normalized right invariant measure on ${\rm B}(2)_\A$.
Note that if we write $b=n(u')t$ where $u'\in\A, t\in{\rm T}(2)_\A$,
then $db=\md t du'$.
We define the local version of
the right invariant measure $db_v$ on ${\rm B}(2)_{k_v}$ similarly.
Then $d_\pr b=\Delta_k^{1/2}\gC_k^{2}db$.
The group $\gl(2)_\A$ has the decomposition
$\gl(2)_\A=\cK(2){\rm B}(2)_\A$.
We choose an invariant measure on $\gl(2)_\A$ by
$dg=d\kappa db$ for $g=\kappa b$.
We define an invariant measure $dg_v$ on $\gl(2)_{k_v}$ similarly.
Then $d_\pr g=\Delta_k^{1/2}\gC_k^{2}dg$.

Recall that we defined $H=\gl(1)\times\gl(2)$.
We regard $\gl(1)$ and $\gl(2)$ as subgroups of $H$.
We choose $dh=\md t dg$ as the normalized invariant measure
on $H_\A=\ma\times\gl(2)_\A$.
We express elements of $B$ as
\begin{equation*}
a(t,p)=\left(t,\twtw 100p\right),\quad
n(u)=\left(1,\twtw 10u1\right).
\end{equation*}
We express element $b\in B_\A$ as
$b=n(u)a(t,p)$ where $t,p\in\ma$ and $u\in\A$.
We choose $db=\md t\md p du$
as a normalized right invariant measure on $B_\A$.
The local measures $dh_v, db_v$ on $H_{k_v}, B_{k_v}$ are defined
in the similar manner and we have
$d_\pr h=\Delta_k^{1/2}\gC_k^{3}dh$
and $d_\pr b=\Delta_k^{1/2}\gC_k^{2}db$.
Let $\gl(2)_\A^1=\{g\in\gl(2)_\A\mid|\det(g)|=1\}$.
We choose a measure $dg^1$ on $\gl(2)_\A^1$ such that
\begin{equation*}
\int_{\gl(2)_\A}f(g)\,dg
=\int_{\R_+}\int_{\gl(2)_\A^1}
	f(\diag(\lamb^{1/2},\lamb^{1/2})g^1)\,dg^1\md\lam.
\end{equation*}
It is well known that the volume of $\gl(2)_\A^1/\gl(2)_k$
with respect to $dg^1$ is $\Delta_k\zeta_k(2)/\gC_k$.

Let $\ac{\cdot}$ is a non-trivial additive character on $\A/k$
and suppose $\ac{x}=\prod_{v\in\gM}\ac{x_v}_v$ for $x=(x_v)\in\A$.
We can choose $\ac{\cdot}$ so that
$\ac{x_v}_v=e^{2\pi ix_v}$ if $v\in\gM_\R$ and
$\ac{x_v}_v=e^{2\pi i(x_v+\overline{x_v})}$ if $v\in\gM_\C$
and fix such a character for all.
We recall the Iwasawa-Tate zeta function.
For $\Phi\in\sS(\A)$, $s\in\C$ and $\omega\in\Omega^1$,
we define
\begin{equation}\label{IT}
\Sigma(\Phi,s,\omega)=\int_{\ma}|t|^s\omega(t)\Phi(t)\md t.
\end{equation}
Let $\Sigma^+(\Phi,s,\omega)$
be the integral obtained from $\Sigma(\Phi,s,\omega)$
by restricting the domain of integration to
$\{t\in\ma \mid |t|\geq1\}$.
The integral $\Sigma^+(\Phi,s,\omega)$
is an entire function of $s$.
We define the Fourier transform $\Phi^\ast$ of $\Phi$ via
$\Phi^\ast(x)=\int_\A\Phi(y)\ac{xy}dy$.
Then by the Poisson summation formula,
we have
\begin{equation*}
\Sigma(\Phi,s,\omega)
=\Sigma^+(\Phi,s,\omega)+\Sigma^+(\Phi^\ast,1-s,\omega^{-1})
+\delta(\omega)\left(\frac{\Phi^\ast(0)}{s-1}-\frac{\Phi(0)}{s}\right).
\end{equation*}
We put
\begin{equation}\label{eq:Sigmazero}
\Sigma_{(0)}(\Phi,\omega)
=\Sigma^+(\Phi,0,\omega)+\Sigma^+(\Phi^\ast,1,\omega^{-1})
	-\delta(\omega)\Phi^\ast(0),
\end{equation}
which is the constant term of
the Laurent expansion of $\Sigma(\Phi,s,\omega)$ at $s=0$.

\section{Global theory for the spaces of binary quadratic forms}
\label{sec:global}

In this section, we define the global zeta functions
for the spaces of binary quadratic forms 
described in Section \ref{sec:pv} and give
some principal part formulae as well as analytic continuations.
These cases are typical example of prehomogeneous vector spaces,
and after Shintani's pioneering work \cite[Chapter1]{shintanib}
the zeta functions were investigated by many authors including
\cite{fsatoa}, \cite{hsaitod}.
Since most of them are written in the classical language,
we reconsidered it in the adelic settings
along the line with F. Sato's formulation \cite{fsatoa}.
Because of Proposition \ref{prop:reducibleparameterization},
our main interest is  the zeta function associated to $(B,W)$.
However as we will see in Lemma \ref{UW}
(which is due to F. Sato),
the zeta functions for $(H,U)$ and $(B,W)$ are
essentially the same, and to give an analytic continuation
of the zeta function it seems to be more natural
to consider representations of reductive groups.
Hence we mainly treat $(H,U)$ in this section.

\subsection{The zeta functions}
We start with the definition of the zeta function for $(H,U)$.
\begin{defn}\label{defn:zetaintHU}
For $\phiu\in\sS(U_\A)$, $s_1,s_2\in\C$ and $\omega_1,\omega_2\in\Omega^1$,
we define
\begin{align*}
X(\phiu,s_1,s_2,\omega_1,\omega_2)
&=	\int_{H_\A/H_k}
		|t|^{s_1}|t\det g|^{2s_2}\omega_1(t)\omega_2(t\det g)
			\sum_{\ti y\in U_k^\sst}\phiu(h\ti y)dh,\\
X^{(1+)}(\phiu,s_1,s_2,\omega_1,\omega_2)
&=	\int_{\underset{|t\det g|\geq1}{H_\A/H_k}}
		|t|^{s_1}|t\det g|^{2s_2}\omega_1(t)\omega_2(t\det g)
			\sum_{\ti y\in U_k^\sst}\phiu(h\ti y)dh,\\
X^{(2+)}(\phiu,s_1,s_2,\omega_1,\omega_2)
&=	\int_{\underset{|\det g|\leq1}{H_\A/H_k}}
		|t|^{s_1}|t\det g|^{2s_2}\omega_1(t)\omega_2(t\det g)
			\sum_{\ti y\in U_k^\sst}\phiu(h\ti y)dh.
\end{align*}
\end{defn}

%
%
For the rest of this section
we consider the analytic continuation of
$X(\phiu,s_1,s_2,\omega_1,\omega_2)$.
When discussing the convergence of integrals
we use the phrase {\em normally convergent}
to mean absolutely and locally uniformly convergent.
The convergence of zeta functions of prehomogeneous vectors space
was investigated by many mathematicians and finally achieved
by H. Saito \cite{hsaitoc}.
For the convergence of the integrals above,
the following lemma holds.
In fact, (1) is contained in \cite{hsaitoc} and
(2) and (3) immediately follow from (1).
\begin{lem}\label{lem:saitoconv}
There exist
$\delta_1>0,\delta_2>0$ such that the following hold.
\begin{enumerate}[{\rm (1)}]
\item\label{saitoconv1}
The integral $X(\phiu,s_1,s_2,\omega_1,\omega_2)$
is normally convergent for $\Re(s_1)>\delta_1$, $\Re(s_2)>\delta_2$.
\item\label{saitoconv2}
The integral $X^{(1+)}(\phiu,s_1,s_2,\omega_1,\omega_2)$
is normally convergent for $\Re(s_1)>\delta_1$.
\item\label{saitoconv3}
The integral $X^{(2+)}(\phiu,s_1,s_2,\omega_1,\omega_2)$
is normally convergent for
$\Re(s_1+2s_2)>\delta_1+2\delta_2$, $\Re(s_2)>\delta_2$.
\end{enumerate}
\end{lem}
\begin{rem}\label{rem:conv}
H. Saito's result is stronger than the above that
the optimum convergence domain is obtained by taking $\delta_1=\delta_2=1$.
In this section we give an alternative proof of this fact.
\end{rem}
Before starting the analysis,
we make two natural assumptions on $\phiu$ for practical purposes.
The first one is:
\begin{asmp}
The test function $\phiu\in\sS(U_\A)$
is of the form $\phiu=\phiw\otimes\phis$,
where $\phiw\in\sS(W_\A)$ and $\phis\in\sS(S_\A)$.
\end{asmp}
Since $W$ and $S$ are both $H$-invariant subspaces,
this is enough for most of the applications.
Let $\cK$ be the standard maximal compact subgroup of $H_\A$ i.e.,
$\cK=\A^0\times\cK(2)$. Let $d\kappa=\md t^0d\kappa_2$
be the measure on $\cK$, so that the total volume of $\cK$ is $1$.
For $\phiu\in\sS(U_\A)$, we define
$\cM_{(\omega_1,\omega_2)}\phiu\in\sS(U_\A)$
by
\begin{equation*}
\cM_{(\omega_1,\omega_2)}\phiu(x)
=	\int_{\cK}
		\omega_1(\chi_1(\kappa))\omega_2(\chi_2(\kappa))\phiu(\kappa x)
	d\kappa.
\end{equation*}
Then we have
$X(\phiu,s_1,s_2,\omega_1,\omega_2)
=X(\cM_{(\omega_1,\omega_2)}\phiu,s_1,s_2,\omega_1,\omega_2)$
and
$\cM_{(\omega_1,\omega_2)}(\cM_{(\omega_1,\omega_2)}\phiu)=
\cM_{(\omega_1,\omega_2)}\phiu$.
Hence, without loss of generality, we will assume:
\begin{asmp}
The Schwartz-Bruhat function $\phiu$ satisfies
$\cM_{(\omega_1,\omega_2)}\phiu=\phiu$.
\end{asmp}
This assumption holds, for example,
if $\phiu$ is $\cK$-invariant and
$\omega_1,\omega_2$ are trivial on $\A^0$
(cf. \cite{kawr}).
This assumption yields the similar assumptions
for the components $\phiw$ and $\phis$.

We now define the zeta function for $(B,W)$.
\begin{defn}\label{zetaintBW}
For $\phiw\in\sS(W_\A)$, $s_1,s_2\in\C$ and $\omega_1,\omega_2\in\Omega^1$,
we define
\begin{align*}
Y(\phiw,s_1,s_2,\omega_1,\omega_2)
&=	\int_{B_\A/B_k}
		|t|^{s_1}|tp|^{2s_2}\omega_1(t)\omega_2(tp)
			\sum_{y\in W_k^\sst}\phiw(by)db,\\
Y^+(\phiw,s_1,s_2,\omega_1,\omega_2)
&=	\int_{\underset{|tp|\geq1}{B_\A/B_k}}
		|t|^{s_1}|tp|^{2s_2}\omega_1(t)\omega_2(tp)
			\sum_{y\in W_k^\sst}\phiw(by)db.
\end{align*}
\end{defn}
For the zeta functions $X(\phiu,s_1,s_2,\omega_1,\omega_2)$
and $Y(\phiw,s_1,s_2,\omega_1,\omega_2)$, the following lemma holds.
Recall that the distribution $\Sigma$ is 
the Iwasawa-Tate zeta function we defined in \eqref{IT}.
(The local version of this lemma
is also given in Lemma \ref{lem:zetadecompbq}.)
\begin{lem}\label{UW}
Let us define $\cR_1\phis\in\sS(\A)$ by $\cR_1\phis(x)=\phis(x,0)$.
We have
\begin{align*}
X(\phiu,s_1,s_2,\omega_1,\omega_2)
&=Y(\phiw,s_1,s_2,\omega_1,\omega_2)\Sigma(\cR_1\phis,2s_1,\omega_1^2),\\
X^{(1+)}(\phiu,s_1,s_2,\omega_1,\omega_2)
&=Y^+(\phiw,s_1,s_2,\omega_1,\omega_2)\Sigma(\cR_1\phis,2s_1,\omega_1^2).
\end{align*}
Especially, the integral $Y(\phiw,s_1,s_2,\omega_1,\omega_2)$
is normally convergent in the region
$\Re(s_1)>\delta_1$, $\Re(s_2)>\delta_2$
and the integral $Y^+(\phiw,s_1,s_2,\omega_1,\omega_2)$
is normally convergent in the region
$\Re(s_1)>\delta_1$.
\end{lem}
\begin{proof}
We consider the first equation.
The second one is proved exactly the same way.
Let
\begin{equation*}
Z'=\{\ti y=(y,\y)\in U\mid Q_1(y)Q_2(y)\y_{21}\not=0, \y_{22}=0\}
\end{equation*}
and $B'=\gl(1)\times\br(2)$.
Then it is easy to see that $U_k^\sst=H_k\times_{B'_k}Z'_k$.
We choose the measure $db'$ on $B_\A'$ as
the product measure of $\md t$ on $\A^\times$ and
$db_2$ on $\br(2)_\A$.
By denoting
$\ti\omega_1=|\cdot|^{s_1}\omega_1$
and $\ti\omega_2=|\cdot|^{2s_2}\omega_2$
we have
\begin{align*}
X(\phiu,s_1,s_2,\omega_1,\omega_2)
&=	\int_{H_\A/B'_{k}}
		\ti\omega_1(t)\ti\omega_2(t\det g)
			\sum_{x\in Z'_k}\phiu(hx)\,dh\\
&=	\int_{B'_{\A}/B'_{k}}
		\ti\omega_1(t)\ti\omega_2(t\det b_2)
			\sum_{x\in Z'_k}\phiu(b'x)\,db'\\
\begin{split}
=	\int_{(\ma/\mk)^3\times\A/k}
		\ti\omega_1(t)\ti\omega_2(tt_1t_2)
		\sum_{y\in W_k^\sst}\phiw((t,\diag(t_1,t_2)n(u))y)\\
\times		\sum_{\y_1\in k^\times}R_1\phis(t_1^{-1}\y_1)\,
			|\frac{t_2}{t_1}|\md t\md t_1\md t_2 du,
\end{split}
\intertext{by changing $t_1$ to $t_1^{-1}$ and after that
$t$ to $tt_1^2$ and $t_2$ to $t_2/t_1$, we have}
&=	\int_{(\ma/\mk)^2\times\A/k}
		\ti\omega_1(t)\ti\omega_2(tt_2)
		\sum_{y\in W_k^\sst}\phiw((a(t,t_2)n(u))y)
			|t_2|\md t\md t_2 du\\
&\quad\times
	\int_{\ma/\mk}\ti\omega_1(t_1^2)
	\sum_{\y_1\in k^\times}R_1\phis(t_1\y_1)\, \md t_1.
\end{align*}
Hence we have the lemma.
\end{proof}

\subsection{Principal part formula I}
In this subsection, we give analytic continuations
of functions $X(\phiu,s_1,s_2,\omega_1,\omega_2)$ and
$Y(\phiw,s_1,s_2,\omega_1,\omega_2)$
to the region $\Re(s_1)>\delta_1$
and find the principal parts in this region.
We first define a singular distribution,
which arises as a principal part of
$Y(\phiw,s_1,s_2,\omega_1,\omega_2)$.
\begin{defn}
For $\phiw\in\sS(W_\A)$, $s\in\C$ and
$\omega\in\Omega^1$,
we put
\begin{equation*}
\Lambda(\phiw,s,\omega)
	=\int_{\ma\times\A}|t|^s\omega(t)\phiw(t,2tu,tu^2)\md t du.
\end{equation*}
\end{defn}
\begin{lem}
The integral $\Lambda(\phiw,s,\omega)$ is normally convergent
for $\Re(s)>1$.
\end{lem}
\begin{proof}
We define the local version of $\Lambda$  by
\begin{equation*}
\Lambda_v(\phiw_{v},s,\omega_v)
	=\int_{k_v^\times\times k_v}
		|t_v|_v^s\omega_v(t_v)\phiw_{v}(t_v,2t_vu_v,t_vu_v^2)\md t_v du_v.
\end{equation*}
For $v\in\gMf$ we let $\phiw_{v,0}$
the characteristic function  of $W_{\co_v}$.
Then by computation we have
\begin{equation*}
\Lambda_v(\phiw_{v,0},s,\omega_v)=
\delta_v(\omega_v)\frac{1-q_v^{-2s}}{(1-q_v^{-s})(1-q_v^{-2s+1})}.
\end{equation*}
By considering the Euler product, we have the lemma.
\end{proof}
We now consider $(B,W)$.
We define the symmetric bilinear forms on $W$ by
%
$[y,z]_W=y_1z_3-2^{-1}y_2z_2+y_3z_1$.
%
Let $\iota$ be the involution on $B$ given by
\begin{equation*}
\iota\colon B\longrightarrow B,
\quad
b=\left(t,\twtw 10up\right)\longmapsto 
b^\iota=\left(t^{-1}p^{-2},\twtw 10up\right).
\end{equation*}
Then we have $[by,b^\iota z]_W=[y,z]_W$
for all $y,z\in W$ and $b\in B$.
\begin{defn}
For $\phiw\in\sS(W_\A)$, we define its Fourier transform 
$\widehat\phiw\in\sS(W_\A)$ by
\begin{equation*}
\widehat{\phiw}(y)=\int_{W_\A}\phiw(z)\ac{[y,z]_W}dz.
\end{equation*}
\end{defn}
For $b\in B_\A$ and $\phiw\in\sS(W_\A)$,
we define $\phiw_{b}\in\sS(W_\A)$ by
$\phiw_{b}(y)=\phiw(by)$.
Then we have
$\widehat{\phiw_{b}}
	=|tp|^{-3}(\widehat\phiw)_{b^\iota}
$.
Finally we define a operator $\cR$ as follows:
\begin{defn}\label{defn:cR}
For $\phiw\in\sS(W_\A)$ we define $\cR\phiw\in\sS(\A)$ by
$\cR\phiw(x)=\int_{\A^2}\phiw(x,u_2,u_3)du_2du_3$.
\end{defn}
The following proposition is an adelic version of
\cite[Lemma 4]{shintanib}.
\begin{prop}\label{pp1}
We have
\begin{align*}
Y(\phiw,s_1,s_2,\omega_1,\omega_2)
&=
	Y^+(\phiw,s_1,s_2,\omega_1,\omega_2)
+	Y^+(\widehat\phiw,s_1,3/2-s_1-s_2,\omega_1,\omega_1^{-2}\omega_2^{-1})\\
&\qquad
+	\frac{\delta(\omega_1^2\omega_2)}{2s_1+2s_2-3}
		\Lambda(\widehat\phiw,s_1,\omega_1)
-	\frac{\delta(\omega_2)}{2s_2}\Lambda(\phiw,s_1,\omega_1)\\
&\qquad
+	\frac{\delta(\omega_2)}{2s_2-2}\Sigma(\cR\phiw,s_1,\omega_1)
-	\frac{\delta(\omega_1^2\omega_2)}{2s_1+2s_2-1}
		\Sigma(\cR\widehat\phiw,s_1,\omega_1).
\end{align*}
\end{prop}
\begin{proof}
Let $\ti\omega_1=|\cdot|^{s_1}\omega_1$
and $\ti\omega_2=|\cdot|^{2s_2}\omega_2$.
We set $Z=\{y\in W\mid y_1(y_2^2-4y_1y_3)=0\}$ and put
\begin{equation*}
I(\phiw,\ti\omega_1,\ti\omega_2)
	=\int_{\underset{|tp|\leq1}{B_\A/B_k}}
		\ti\omega_1(t)\ti\omega_2(tp)
			\left(
			\sum_{y\in Z_k}\widehat{\phiw_{b}}(y)
			-\sum_{y\in Z_k}\phiw_{b}(y)
		\right)
	db.
\end{equation*}
Then since $W_k=W_k^\sst\amalg Z_k$,
by the Poisson summation formula we have
\begin{multline*}
Y(\phiw,s_1,s_2,\omega_1,\omega_2)
=
	Y^+(\phiw,s_1,s_2,\omega_1,\omega_2)\\
+	Y^+(\widehat\phiw,s_1,3/2-s_1-s_2,\omega_1,\omega_1^{-2}\omega_2^{-1})
+I(\phiw,\ti\omega_1,\ti\omega_2).
\end{multline*}
We consider $I(\phiw,\ti\omega_1,\ti\omega_2)$.
Let
\begin{equation*}
Z_1=\{y\in W\mid y_{1}\not=0,y_{2}^2-4y_{1}y_{3}=0\},\quad
Z_2=\{y\in W\mid y_{1}=0\}.
\end{equation*}
Obviously we have $Z_k=Z_{1k}\amalg Z_{2k}$.
We put $B_\A^{-}=\{b\in B_\A\mid |tp|\leq1\}$
and define
\begin{align*}
I_1(\phiw,\ti\omega_1,\ti\omega_2)
&	=\int_{{B^{-}_\A/B_k}}
		\ti\omega_1(t)\ti\omega_2(tp)|tp|^{-3}
			\sum_{y\in Z_{1k}}\widehat\phiw(b^\iota y)
	db,\\
I_2(\phiw,\ti\omega_1,\ti\omega_2)
&	=\int_{{B^{-}_\A/B_k}}
		\ti\omega_1(t)\ti\omega_2(tp)
			\sum_{y\in Z_{1k}}\phiw(by)
	db,\\
I_3(\phiw,\ti\omega_1,\ti\omega_2)
&	=\int_{{B^{-}_\A/B_k}}
		\ti\omega_1(t)\ti\omega_2(tp)
			\left(
			\sum_{y\in Z_{2k}}\widehat{\phiw_{b}}(y)
			-\sum_{y\in Z_{2k}}\phiw_{b}(y)
		\right)
	db.
\end{align*}
Then $I(\phiw,\ti\omega_1,\ti\omega_2)=
(I_1-I_2+I_3)(\phiw,\ti\omega_1,\ti\omega_2)$.
Let us first consider $I_1$.
%
%
We let $w=(1,0,0)\in Z_{1k}$
and $B_w=\{a(1,p)\mid p\in\mathbb G_m\}$,
which is the stabilizer of $w$ in $B$.
It is easy to see that $Z_{1k}=B_kw$.
Hence we have
\begin{align*}
I_1(\phiw,\ti\omega_1,\ti\omega_2)
&	=\int_{{B^{-}_\A/B_{wk}}}
		\ti\omega_1(t)\ti\omega_2(tp)|tp|^{-3}
			\widehat\phiw(b^\iota w)
	db\\
&	=\int_{\underset{|tp|\leq1}{(\ma\times\ma/\mk)}\times\A}
		\ti\omega_1(t)\ti\omega_2(tp)|tp|^{-3}
			\widehat\phiw(1/tp^2,2u/tp^2,u^2/tp^2)
	\, \md t\md p du
\intertext{by changing $t$ to $t^{-1}p^{-2}$,}
&	=\int_{\underset{|tp|\geq1}{(\ma\times\ma/\mk)}\times A}
		(\ti\omega_1\ti\omega_2)(t^{-1})|t|^3
			\widehat\phiw(t,2tu,tu^2)
	\cdot(\ti\omega_1^2\ti\omega_2)(p^{-1})|p|^3	\, \md t\md p du.
\end{align*}
Since
\begin{equation*}
\int_{\underset{|tp|\geq1}{\ma/\mk}}
	(\ti\omega_1^2\ti\omega_2)(p^{-1})|p|^3	\, \md p
=\frac{\delta(\omega_1^2\omega_2)}{2s_1+2s_2-3}(\ti\omega_1^2\ti\omega_2)(t)|t|^{-3},
\end{equation*}
we have
$I_1(\phiw,\ti\omega_1,\ti\omega_2)=(2s_1+2s_2-3)^{-1}
\delta(\omega_1^2\omega_2)\Lambda(\widehat\phiw,s,\omega_1)$.
Similarly we have
$I_2(\phiw,\ti\omega_1,\ti\omega_2)
=(2s_2)^{-1}{\delta(\omega_2)}\Lambda(\phiw,s_1,\omega_1)$.
We finally consider $I_3(\phiw,\ti\omega_1,\ti\omega_2)$.
We have
\begin{equation*}
I_3(\phiw,\ti\omega_1,\ti\omega_2)
=	\int_{T^{-}_\A/T_k}
	\ti\omega_1(t)\ti\omega_2(tp)|p|
J(\phiw_{a(t,p)})\md t
\end{equation*}
where we put $T_A^{-}=B_\A^-\cap T_\A$ and
\begin{equation*}
J(\phiw)=\int_{N_\A/N_k}
\left(			\sum_{y\in S_{2k}}\widehat\phiw\bigl(n(u)y\bigr)
				-\sum_{y\in S_{2k}}\phiw\bigl(n(u)y\bigr)
\right)					du.
\end{equation*}
We consider $J(\phiw)$.
We define $\cR_3\phiw\in\sS(\A)$ by
$\cR_3\phiw(x)=\phiw(0,0,x)$.
By dividing the index set $Z_{2k}$ into
\begin{equation*}
Z_{2k}=
	\{(0,y_{2},y_{3})\mid y_{2}\in\mk, y_{3}\in k\}
	\amalg
	\{(0,0,y_{3})\mid y_{3}\in k\},
\end{equation*}
and performing integration separately, we have
\begin{equation*}
J(\phiw)
=	\sum_{\alpha\in k}\cR_3\widehat\phiw(\alpha)
-	\sum_{\alpha\in k}\cR_3\phiw(\alpha)
+	\sum_{\alpha\in\mk}\int_{\A}\widehat\phiw(0,\alpha,u)du
-	\sum_{\alpha\in\mk}\int_{\A}\phiw(0,\alpha,u)du.
\end{equation*}
Since the equality
$
\sum_{\alpha\in k}\int_{\A}\widehat\phiw(0,\alpha,u)du=
\sum_{\alpha\in k}\int_{\A}\phiw(0,\alpha,u)du$
holds by the \psf, we have
\begin{equation*}
J(\phiw)
=	\left(\sum_{\alpha\in k}\cR_3\widehat\phiw(\alpha)
		-\int_{\A}\cR_3\widehat\phiw(u)du\right)
-	\left(\sum_{\alpha\in k}\cR_3\phiw(\alpha)
		-\int_{\A}\cR_3\phiw(u)du\right).
\end{equation*}
Again by using the \psf, we have
\begin{equation*}
J(\phiw)
=	\sum_{\alpha\in\mk}\cR\phiw(\alpha)-\sum_{\alpha\in\mk}\cR\widehat\phiw(\alpha).
\end{equation*}
Note that we define the operator $\cR$ in Definition \ref{defn:cR}.
Now we can easily see
\begin{equation*}
J(\phiw_{a(t,p)})
=	|t^{-2}p^{-3}|\sum_{\alpha\in\mk}\cR\phiw(t\alpha)
	-|t^{-1}p^{-2}|\sum_{\alpha\in\mk}\cR\widehat\phiw(t^{-1}p^{-2}\alpha)
\end{equation*}
and hence,
\begin{align*}
I_3(\phiw,\ti\omega_1,\ti\omega_2)
&=	\int_{T^{-}_\A/T_k}
	\ti\omega_1(t)\ti\omega_2(tp)|tp|^{-2}\sum_{\alpha\in \mk}\cR\phiw(t\alpha)
		\, \md t\md p\\
&\quad
-	\int_{T^{-}_\A/T_k}
	\ti\omega_1(t)\ti\omega_2(tp)|tp|^{-1}\sum_{\alpha\in \mk}
	\cR\widehat\phiw(t^{-1}p^{-2}\alpha)
		\, \md t\md p.
\end{align*}
Now straightforward calculations show that
the integrals above equal to
\begin{equation*}
(2s_2-2)^{-1} \delta(\omega_2) \Sigma(\cR\phiw,s_1,\omega_1)
\quad\text{and}\quad
(2s_1+2s_2-1)^{-1}\delta(\omega_1^2\omega_2)\Sigma(\cR\widehat\phiw,s_1,\omega_1),
\end{equation*}
respectively.
These give the desired description.
\end{proof}

As a corollary to this proposition, we obtain the following.
\begin{cor}\label{ac1}
The function
$s_2(s_2-1)(2s_1+2s_2-1)(2s_1+2s_2-3)X(\phiu,s_1,\omega_1,s_2,\omega_2)$
is holomorphically continued to the region $\Re(s_1)>\delta_1$.
Also the following functional equation holds:
\begin{equation*}
X(\phiw\otimes\phis,s_1,s_2,\omega_1,\omega_2)
=X(\widehat{\phiw}\otimes\phis,s_1,3/2-s_1-s_2,\omega_1,\omega_1^{-2}\omega_2^{-1}).
\end{equation*}
\end{cor}
\subsection{Principal part formula II}
In this subsection, we give analytic continuation
of the function $X(\phiu,s_1,s_2,\omega_1,\omega_2)$ to the region
$\Re(s_1+2s_2)>\delta_1+2\delta_2, \Re(s_2)>\delta_2$
and find the principal parts in this region.
\begin{defn}
Let $U_{(1)}=\{\ti y\in U_k^\sst\mid k(\ti y)\cong k\times k\}$
and $U_{(2)}=U_k^\sst\setminus U_{(1)}$.
For $i=1,2$, we define
\begin{align*}
X_{(i)}(\phiu,s_1,s_2,\omega_1,\omega_2)
&=	\int_{H_\A/H_k}
	|t|^{s_1}|t\det g|^{2s_2}
		\omega_1(t)\omega_2(t\det g)
			\sum_{x\in U_{(i)}}\phiu(hx)dh,\\
X^{(2+)}_{(i)}(\phiu,s_1,s_2,\omega_1,\omega_2)
&=	\int_{\underset{|\det g|\leq1}{H_\A/H_k}}
	|t|^{s_1}|t\det g|^{2s_2}
		\omega_1(t)\omega_2(t\det g)
			\sum_{x\in U_{(i)}}\phiu(hx)dh.
\end{align*}
\end{defn}
Since $U_k^\sst=U_{(1)}\amalg U_{(2)}$
we have $X=X_{(1)}+X_{(2)}$ and
$X^{(2+)}=X^{(2+)}_{(1)}+X^{(2+)}_{(2)}$.
By Lemma \ref{lem:saitoconv} (\ref{saitoconv3}),
$X^{(2+)}_{(i)}(\phiu,s_1,s_2,\omega_1,\omega_2)$ is holomorphic for
$\Re(s_1+2s_2)>\delta_1+2\delta_2,\Re(s_2)>\delta_2$.
We consider $X_{(1)}$ and $X_{(2)}$ separately.

We define the bilinear forms on $S$ by
$[\y,\z]_S=\y_1\z_2-\y_2\z_1$.
Let $\iota$ be the involution on $H$ given by
$(t,g)^\iota=(t,-(\det g)^{-1}g)$.
Then we have $[h\y,h^\iota \z]_S=[\y,\z]_S$
for all $\y,\z\in S$ and $h\in H$.
\begin{defn}
For $\phis\in\sS(S_{\A})$, we define its Fourier transform 
$\widehat\phis$ by
\begin{equation*}
\widehat{\phis}(\y)=\int_{S_{\A}}\phis(\z)\ac{[\y,\z]_S}d\z.
\end{equation*}
\end{defn}
%
%

We first consider $X_{(2)}(\phiu,s_1,s_2,\omega_1,\omega_2)$.
Since the zeta function associated with  $(H,W)$
appears as a term of the residue at $s_1=1$
(this observation is due to Shintani \cite{shintanib}),
we recall the definition of the zeta function for the space.
In this case, due to the convergence problem,
we are not allowed to let the index set of the summation
in the integral to be $W_{k}\tilde{}$.
We use $W_{(2)}$ defined below instead.
Let $dh^1=\md t dg^1$ for $h^1=(t,g^1), t\in\ma, g^1\in\gl(2)_\A^1$.
\begin{defn}
Let $W_{(2)}=\{y\in W_k\tilde{}\mid k(y)\not\cong k\times k\}$.
For $\phiw\in\sS(W_\A)$, $s_1,s_2\in\C$ and $\omega_1,\omega_2\in\Omega^1$,
we define
\begin{equation*}
U(\phiw,s_1,s_2,\omega_1,\omega_2)
=\int_{H_\A^1/H_k}
|t|^{s_1+2s_2}\omega_1(t)\omega_2(t\det g^1)
	\sum_{y\in W_{(2)}}\phiw(h^1y)dh^1.
\end{equation*}
\end{defn}
The function $U(\phiw,s_1,s_2,\omega_1,\omega_2)$
is essentially in one variable
but we use the above definition for conveniences.
It is well known that the integral converges if $\Re(s_1+2s_2)>\delta_3$
for some $\delta_3$. By changing $\delta_1$ and $\delta_2$ if necessary
we may assume $\delta_3=\delta_1+2\delta_2$.
For $p\in\mathbb G_m$, we have $(tp^2,p^{-1}g)y=(t,g)y$.
Hence by computation we could see:
\begin{lem}
\begin{equation*}
U(\phiw,s_1,s_2,\omega_1,\omega_2)
=U(\phiw,1-s_1,s_1+s_2-1/2,\omega_1^{-1},\omega_1^2\omega_2).
\end{equation*}
\end{lem}
\begin{prop}\label{pp2}
\begin{align*}
&X_{(2)}(\phiw\otimes\phis,s_1,s_2,\omega_1,\omega_2)\\
&\quad=X_{(2)}^{(2+)}(\phiw\otimes\phis,s_1,s_2,\omega_1,\omega_2)
+X_{(2)}^{(2+)}(\phiw\otimes\widehat{\phis},1-s_1,s_1+s_2-1/2,\omega_1^{-1},\omega_1^2\omega_2)\\
&\quad
+\frac{\widehat{\phis}(0)}{s_1-1}
	U(\phiw,s_1,s_2,\omega_1,\omega_2)
-\frac{\phis(0)}{s_1}
	U(\phiw,s_1,s_2,\omega_1,\omega_2).
\end{align*}
\end{prop}
\begin{proof}
Let $\ti\omega_1=|\cdot|^{s_1}\omega_1$
and $\ti\omega_2=|\cdot|^{2s_2}\omega_2$.
Since $y(v)$ is a irreducible quadratic polynomial
for $\ti y=(y,\y)\in U_{(2)}$, $R_1(\ti y)=0$ implies $\y=0$,
i.e. , $U_{(2)}=W_{(2)}\times\{\y\in S_{k}\mid \y\not=(0,0)\}$.
Hence by the Poisson summation formula,
\begin{align*}
&X_{(2)}(\phiu,s_1,s_2,\omega_1,\omega_2)
-X_{(2)}^{(2+)}(\phiu,s_1,s_2,\omega_1,\omega_2)\\
&\qquad
=\int_{\underset{|\det g|\geq1}{H_\A/H_k}}
		\ti\omega_1(t)\ti\omega_2(t\det g)|\det g|
		\sum_{\ti y\in U_{(2)}}
			\phiw(hy)\phis(h^\iota\y)dh\\
&\qquad\quad
+\int_{\underset{|\det g|\geq1}{H_\A/H_k}}
		\ti\omega_1(t)\ti\omega_2(t\det g)
			\sum_{y\in W_{(2)}}\phiw(hy)
		\left(|\det g|\widehat{\phis}(0)-\phis(0)\right)dh.
\end{align*}
It is easy to see that the second integral in the formula above
corresponds to the last second terms in the formula of the proposition.
We consider the first integral in the formula, which is equal to
\begin{equation*}
\int_{\underset{|\det(g)|\leq1}{H_\A/H_k}}
		\omega_1(t)\omega_2(t/\det g)|\det g|^{-1}
		\sum_{\ti y\in U_{(2)}}
			\phiw(h^\iota y)\phis(h\y)dh.
\end{equation*}
For $y\in W$ we have
%
$h^\iota y
=(t,-(\det g)^{-1}g)y=(t(\det g)^{-2},g)y$.
%
With this relation in mind, by changing $t$ to $t(\deg g)^2$
we could see that the integral is equal to
$X_{(2)}^{(2+)}(\phiw\otimes\widehat{\phis},1-s_1,s_1+s_2-1/2,\omega_1^{-1},\omega_1^2\omega_2)$.
\end{proof}

We next consider $X_{(1)}^{(2+)}(\phiu,s_1,s_2,\omega_1,\omega_2)$.
We introduce a function which plays an important role
in our analysis.
\begin{defn}
For $\phis\in\sS(S_\A)$, $s\in\C$ and $\omega\in\Omega^1$, let
\begin{align*}
\Sigma_2(\phis,s,\omega)
&=\int_{\ma\times\ma}|t_1t_2|^s\omega(t_1t_2)\phis(t_1,t_2)\md t_1\md t_2,\\
\Sigma_2^+(\phis,s,\omega)
&=\int_{\underset{|t_1t_2|\geq1}{\ma\times\ma}}
	|t_1t_2|^s\omega(t_1t_2)\phis(t_1,t_2)\md t_1\md t_2.
\end{align*}
\end{defn}
Note that the function $\Sigma_2^+(\phis,s,\omega)$ is entire.
We put $Z=\{(x_{21}x_{22})\in S\mid x_{21}x_{22}=0\}$.
Let
\begin{align*}
K(\phis,t_1,t_2)
&=	|t_1t_2|^{-1}\sum_{{x_2\in Z_{k}}}
			\widehat\phis(t_2^{-1}x_{21},t_1^{-1}x_{22})
-		\sum_{{x_2\in Z_{k}}}
			\phis(t_1x_{21},t_2x_{22}),\\
J'(\phis,s,\omega)
&=\int_{\underset{|t_1t_2|\leq1}{(\ma/\mk)^2}}
	|t_1t_2|^s\omega(t_1t_2)K(\phis,t_1,t_2)
\md t_1\md t_2.
\end{align*}
Then by the Poisson summation formula, we have
\begin{equation*}
\Sigma_2(\phis,s,\omega)
=\Sigma_2^+(\phis,s,\omega)+\Sigma_2^+(\widehat\phis,1-s,\omega^{-1})
+J'(\phis,s,\omega).
\end{equation*}
For $\phis\in\sS(S_\A)$,
we define $\cR_1\phis,\cR_2\phis\in\sS(\A)$
by $\cR_1\phis(x)=\phis(x,0), \cR_2\phis(x)=\phis(0,x)$.
Then the following holds.
Note that we define the distribution
$\Sigma_{(0)}$ in \eqref{eq:Sigmazero}.
\begin{lem}\label{Sigma2}
We have $J'(\phis,s,\omega)=\delta(\omega)J''(\phis,s)$ where
\begin{align*}
J''(\phis,s)
&=	\frac{\widehat\phis(0)}{(s-1)^2}
+	\frac{1}{s-1}
	\left(
		\Sigma_{(0)}(\cR_1\widehat\phis,0)+\Sigma_{(0)}(\cR_2\widehat\phis,0)
	\right)\\
&\quad
+\frac{\phis(0)}{s^2}
-	\frac{1}{s}
	\left(
		\Sigma_{(0)}(\cR_1\phis,0)+\Sigma_{(0)}(\cR_2\phis,0)
	\right).
\end{align*}
\end{lem}
\begin{proof}
We shall calculate the integral $J'$
by dividing it to the three integrals
\begin{equation*}
J'_0=\int_{|t_1|\leq1,|t_2|\leq1},
\quad
J'_1=\int_{1\leq |t_1|\leq|t_2|^{-1}}, 
\quad
J'_2=\int_{1\leq |t_2|\leq|t_1|^{-1}}.
\end{equation*}
Then by suitable uses of the Poisson summation formula for each integral,
we have $J'_i(\phis,\omega)=\delta(\omega)J''_i(\phis,s)$
where
\begin{align*}
J''_0(\phis,s)
&=	\frac{\widehat\phis(0)}{(s-1)^2}+\frac{\phis(0)}{s^2}
-	\frac{(\cR_1\phis)^\ast(0)+(\cR_2\phis)^\ast(0)}{s(s-1)}\\
&\qquad
+	\frac{1}{s-1}
	\left(
		\Sigma_+(\cR_1\widehat\phis,1-s)+\Sigma_+(\cR_2\widehat\phis,1-s)
	\right)\\
&\qquad
-	\frac{1}{s}
	\left(
		\Sigma_+((\cR_1\phis)^\ast,1-s)+\Sigma_+((\cR_2\phis)^\ast,1-s)
	\right),\\
J''_1(\phis,s)
&
=	\frac{1}{s-1}
	\left(
		\Sigma_+((\cR_2\widehat\phis)^\ast,1)
		+\Sigma_+(\cR_1\widehat\phis,0)
		-\Sigma_+(\cR_1\widehat\phis,1-s)
	\right)\\
&\qquad
-	\frac{1}{s}
	\left(
		\Sigma_+((\cR_2\phis)^\ast,1)
		+\Sigma_+(\cR_1\phis,0)
		-\Sigma_+((\cR_2\phis)^\ast,1-s)
	\right),\\
J''_2(\phis,s)
&
=	\frac{1}{s-1}
	\left(
		\Sigma_+((\cR_1\widehat\phis)^\ast,1)
		+\Sigma_+(\cR_2\widehat\phis,0)
		-\Sigma_+(\cR_2\widehat\phis,1-s)
	\right)\\
&\qquad
-	\frac{1}{s}
	\left(
		\Sigma_+((\cR_1\phis)^\ast,1)
		+\Sigma_+(\cR_2\phis,0)
		-\Sigma_+((\cR_1\phis)^\ast,1-s)
	\right).
\end{align*}
Since
$(\cR_1\widehat\phis)^\ast(0)=(\cR_2\phis)^\ast(0)$ and
$(\cR_2\widehat\phis)^\ast(0)=(\cR_1\phis)^\ast(0)$,
by adding all them up, we have the formula.
\end{proof}
\begin{defn}
We define $\phis_{n(u)}\in\sS(S_\A)$ by
$\phis_{n(u)}(\y)=\phis(n(u)\y)$.
We put
\begin{align*}
\Pi_1(\phiw,s,\omega)
&=\int_{\ma\times\A}|t|^s\omega(t)\phiw(0,t,tu)du,\\
\Pi_2(\phiu,s,\omega)
&=\int_{\ma\times\A}|t|^s\omega(t)\phiw(0,t,tu)\Sigma_{(0)}(\cR_2(\phis_{n(u)}),0)\md tdu.
\end{align*}%
\end{defn}
It is easy to see that these integrals are normally convergent for $\Re(s)>2$.
\begin{prop}\label{pp3}
We have
\begin{align*}
&X_{(1)}(\phiw\otimes\phis,s_1,s_2,\omega_1,\omega_2)
=X_{(1)}^{(2+)}(\phiw\otimes\phis,s_1,s_2,\omega_1,\omega_2)\\
&\quad+X_{(1)}^{(2+)}(\phiw\otimes\widehat{\phis},1-s_1,s_1+s_2-1/2,\omega_1^{-1},\omega_1^2\omega_2)
+2^{-1}{\delta(\omega_1)}J(\phiu,s_1,s_2,\omega_1\omega_2),
\end{align*}
where $J(\phiu,s_1,s_2,\omega_1\omega_2)$ is given by

\begin{align*}
&	\frac 1{(s_1-1)^2}\widehat{\phis}(0)
		\Pi_1(\phiw,s_1+2s_2,\omega_1\omega_2)
+	\frac 1{s_1^2}\phis(0)\Pi_1(\phiw,s_1+2s_2,\omega_1\omega_2)\\
&
+	\frac 1{s_1-1}\Sigma_{(0)}(\cR_1\widehat{\phis},0)
\Pi_1(\phiw,s_1+2s_2,\omega_1\omega_2)
+	\frac 1{s_1-1}\Pi_2(\phiw\otimes\widehat{\phis},s_1+2s_2,\omega_1\omega_2)\\
&
-	\frac 1{s_1}\Sigma_{(0)}(\cR_1\phis,0)\Pi_1(\phiw,s_1+2s_2,\omega_1\omega_2)
-	\frac 1{s_1}\Pi_2(\phiu,s_1+2s_2,\omega_1\omega_2).\\
\end{align*}
\end{prop}
\begin{proof}
As before we put $\ti\omega_1=|\cdot|^{s_1}\omega_1$
and $\ti\omega_2=|\cdot|^{2s_2}\omega_2$.
Let $U_{(1)}\ni w=(w_1,w_2)$, where $w_1=(0,1,0),w_2=(1,1)$.
We put $B'=\gl(1)\times{\mathrm B}(2)$.
We choose the measure $db'$ on $B_\A'$ as
the product measure of $\md t$ on $\A^\times$ and
$db_2$ on ${\mathrm B}(2)_\A$.
Let $\phiw_{n(u)}(y)=\phiw(n(u)y)$.
Then since $U_{(1)}=H_{k}w$ and ${}^\#(H_{wk})=2$, 
$X_{(1)}(\phiu,s_1,s_2,\omega_1,\omega_2)$ equals to
\begin{align*}
&\frac12
	\int_{H_\A}
		\ti\omega_1(t)\ti\omega_2(t\det g)
			\phiw(hw_1)\phis(hw_2)dh\\
&=\frac12
	\int_{B_\A'}
		\ti\omega_1(t)\ti\omega_2(t\det b_2)
			\phiw(b'w_1)\phis(b'w_2)db'\\
&=\frac12
	\int_{(\ma)^3\times\A}
		\ti\omega_1(t)\ti\omega_2(tt_1t_2)
			\phiw_{n(u)}(0,tt_1t_2,0)
			\phis_{n(u)}(\frac{1}{t_{1}},\frac{1}{t_{2}})
	\md t\md t_1\md t_2 du.
\end{align*}
By changing $t_1, t_2$ to $1/t_1,1/t_2$ and after that
$t$ to $tt_1t_2$, we have
\begin{equation*}
X_{(1)}(\phiu,s_1,s_2,\omega_1,\omega_2)
=\frac12
	\int_{\ma\times \A}\ti\omega_1\ti\omega_2(t_2)
			\phiw_{n(u)}(0,t_2,0)
			\Sigma_2(\phis_{n(u)},s_1,\omega_1)
	\md t_2 du.
\end{equation*}
Also the same modification yields
\begin{equation*}
X_{(1)}^{(2+)}(\phiu,s_1,s_2,\omega_1,\omega_2)
=\frac12
	\int_{\ma\times \A}\ti\omega_1\ti\omega_2(t_2)
			\phiw_{n(u)}(0,t_2,0)
			\Sigma_2^+(\phis_{n(u)},s_1,\omega_1)
	\md t_2 du.
\end{equation*}
Now by applying Lemma \ref{Sigma2},
we obtain the desired formula.
\end{proof}
As a corollary to Propositions \ref{pp2} and \ref{pp3},
we also obtain the following.
\begin{cor}\label{ac2}
The function $s_1^2(s_1-1)^2X(\phiu,s_1,s_2,\omega_1,\omega_2)$
is holomorphically continued to the region
$\Re(s_1+2s_2)>\delta_1+2\delta_2, \Re(s_2)>\delta_2$.
Also the following functional equation holds.
\begin{equation*}
X(\phiw\otimes\phis,s_1,s_2,\omega_1,\omega_2)
=X(\phiw\otimes\widehat{\phis},1-s_1,s_1+s_2-1/2,\omega_1^{-1},\omega_1^2\omega_2).
\end{equation*}
\end{cor}

\begin{rem}
This functional equation is nothing but those of
$L$-functions of quadratic extensions of $k$.
For this fact, see Proposition \ref{prop:orbitalzeta}.
\end{rem}

\subsection{Analytic continuation}
By putting together we have obtained in the two previous subsections,
we obtain meromorphic continuation of the global zeta function.
\begin{thm}\label{thm:fneqHU}
Let 
$\ti X(\phiu,s_1,s_2,\omega_1,\omega_2)
=s_1^2(s_1-1)^2s_2(s_2-1)(2s_1+2s_2-1)(2s_1+2s_2-3)
X(\phiu,s_1,s_2,\omega_1,\omega_2)$.
Then $\ti X(\phiu,s_1,s_2,\omega_1,\omega_2)$
is holomorphically continued to all $\C^2$.
Also the zeta function satisfies the following functional equations:
\begin{multline*}
X(\phiw\otimes\phis,s_1,s_2,\omega_1,\omega_2)
=X(\widehat{\phiw}\otimes\phis,s_1,3/2-s_1-s_2,\omega_1,\omega_1^{-2}\omega_2^{-1})\\
=X(\phiw\otimes\widehat{\phis},1-s_1,s_1+s_2-1/2,\omega_1^{-1},\omega_1^2\omega_2)
=X(\widehat{\phiw}\otimes\widehat{\phis},1-s_1,1-s_2,\omega_1^{-1},\omega_2^{-1}).
\end{multline*}
\end{thm}
\begin{proof}
By Corollaries \ref{ac1} and \ref{ac2}, 
$\ti X(\phiu,s_1,s_2,\omega_1,\omega_2)$
can be continued holomorphically to the tube domain
\begin{equation*}
\cD=
\{(s_1,s_2)\in\C^2\mid \Re(s_1)>\delta_1\}\cup
\{(s_1,s_2)\in\C^2\mid \Re(s_1+2s_2)>\delta_1+2\delta_2,\Re(s_2)>\delta_2\}.
\end{equation*}
Since the convex hull of $\cD$ coincides with $\C^2$,
$\ti X(\phiu,s_1,s_2,\omega_1,\omega_2)$
can be continued holomorphically to the whole $\C^2$
(cf. \cite[Theorem 2.5.10]{hormander}).
The functional equation is obvious.
\end{proof}

As a corollary to this theorem,
combined with the Hartogs theorem \cite[Theorem 2.3.2]{hormander},
we can strength
Lemmata \ref{lem:saitoconv} and \ref{UW} as follows
(see Remark \ref{rem:conv}.)
\begin{cor}\label{cor:conv}
The statements of Lemmata \ref{lem:saitoconv} and \ref{UW}
hold with $\delta_1=\delta_2=1$.
\end{cor}

\section{Archimedean local theory for the spaces of binary quadratic forms}
\label{sec:local}
In this section we describe the gamma factor of the functional equation
of the local zeta functions.
This is used in the proof of Theorem \ref{thm:fneqDHU}
to determine the gamma factor
of the functional equation of the Dirichlet series.
We study some non-archimedean local theory in Appendix \ref{sec:non-arch}.

\subsection{The definition of unramified local zeta functions}
In this subsection $v\in\gM$ is arbitrary.
For convenience we introduce an index set of
the orbits $H_{k_v}\backslash U_{k_v}^\sst$.
\begin{defn}
For each $v\in\gM$ we let $\type_v$ be the index set 
for the set of orbits $H_{k_v}\backslash U_{k_v}^\sst$.
By Proposition \ref{prop:rod2} $\type_v$
corresponds bijectively to the set of isomorphism classes of
separable quadratic algebras of $k_v$.
For $j_v\in\type_v$,
we denote by $U_{k_v,j_v}$ or simply $U_{j_v}$
the $H_{k_v}$-orbit in $U^\sst_{k_v}$ corresponding to $j_v$.
For $v\in\gM_\R$, we further let $\type_v=\{1,2\}$
where the orbit corresponding to $\R\times\R$ is
labeled $1$ and $\C$ labeled $2$.
For $v\in\gM_\C$ we let $\type_v=\{1\}$.
Hence
$U_{\R,1}=\{\ti y\in U_\R^\sst \mid R_2(\ti y)>0\}$,
$U_{\R,2}=\{\ti y\in U_\R^\sst \mid R_2(\ti y)<0\}$
and $U_{\C,1}=U_\C^\sst$.
\end{defn}
The unramified local zeta function is defined as follows:
\begin{defn}
Let $j_v\in\type_v$. Take an arbitrary element $\ti y\in U_{j_v}$.
For $\phiu_v\in\sS(U_{k_v})$ and $s_1,s_2\in\C$, we put
\begin{align*}
\vcX_{v,j_v}(\phiu_v,s_1,s_2)
&=	\int_{H_{k_v}}
		|R_1(h_v\ti y)|_v^{s_1}|R_2(h_v\ti y)|_v^{s_2}
		\phiu_v(h_v\ti y)dh_v,\\
\mho_{v,j_v}(\phiu_v,s_1,s_2)
&=	\int_{U_{j_v}}
		|R_1(\ti y_v)|_v^{s_1-1}|R_2(\ti y_v)|_v^{s_2-1}
		\phiu_v(\ti y_v)d\ti y_v.
\end{align*}
By the uniqueness of the invariant measure,
these coincide up to a positive constant.
We put $\vcX_{v,j_v}=b_{v,j_v}\mho_{v,j_v}$.
\end{defn}
These integrals converges at least $\Re(s_1)>1,\Re(s_2)>1$.
The next lemma follows from a straightforward computation
of the Jacobian determinant of the double cover
$H_{k_v}\ni h\mapsto h\ti y\in U_{j_v}$ for suitable coordinate systems.
\begin{lem}\label{lem:b_v}
Let $v\in\gM_\infty$.
We have $b_{v,j_v}=2\Gamma_{k_v}(2)$.
\footnote{We define $\Gamma_\R(s)$ and $\Gamma_\C(s)$ in
Section \ref{sec:measure}.}
Especially it does not depend on $j_v$.
\end{lem}
We define the local Fourier transform
exactly same way as the global one.
\begin{defn}\label{defn:fourierUlocal}
Let $\phiu_v\in\sS(U_{k_v})$.
We define the Fourier transform of $\phiu_v$ by
\begin{equation*}
\widehat{\phiu_v}(\ti y_v)=\int_{U_{k_v}}\phiu_v(z_v)
\ac{[y_v,z_v]_W+[\y_v,\z_v]_S}_v\, d\ti z_v.
\end{equation*}

\end{defn}
\subsection{Functional equations at infinite places}
In this subsection we give the functional equations
of local zeta functions at $\gM_\infty$.
As in the global situation Section \ref{sec:global},
the local zeta functions satisfy $2$ kinds of functional equation.
However here we deal with only $1$ kind
which we need to prove density theorems.
For other types, see Appendix \ref{sec:remainfneq}.

\begin{prop}\label{prop:localfneq}
We assume $v\in\gM_\infty$.
Let $\ti\Gamma(s_1,s_2)=\Gamma(s_1)^2\Gamma(s_2)\Gamma(s_1+s_2-1/2)$.
\begin{enumerate}[{\rm (1)}]
\item
The function $\mho_{v,j_v}(\phiu_v,s_1,s_2)$
becomes entire after multiplied by $\ti\Gamma(s_1,s_2)^{-1}$.
\item
Let $\phiu\in\sS(U_\R)$. The functional equation for $\R$ is:
\begin{equation*}
\begin{pmatrix}
\mho_{\R,1}(\widehat\phiu,s_1,s_2)\\
\mho_{\R,2}(\widehat\phiu,s_1,s_2)\\
\end{pmatrix}
=\frac{2^{2s_2-s_1}}{\pi^{3s_1+2s_2-1/2}}
\ti\Gamma(s_1,s_2)D_{\R}(s_1,s_2)
\begin{pmatrix}
\mho_{\R,1}(\phiu,1-s_1,1-s_2)\\
\mho_{\R,2}(\phiu,1-s_1,1-s_2)\\
\end{pmatrix}
\end{equation*}
where the $2\times 2$ matrix $D_\R(s_1,s_2)=(d_{\R,i,j}(s_1,s_2))$
is given by
\begin{equation*}
\begin{pmatrix}
2\cos^2(s_1\pi/2)\sin((s_1+2s_2)\pi/2)	&\sin(s_1\pi)\cos(s_1\pi/2)\\
\sin(s_1\pi)\cos(s_1\pi/2)		&\sin(s_1\pi)\cos((s_1+2s_2)\pi/2)\\
\end{pmatrix}.
\end{equation*}
\item
Let $\phiu\in\sS(U_\C)$. The functional equation for $\C$ is:
\begin{equation*}
\mho_{\C,1}(\widehat\phiu,s_1,s_2)
=\frac{2^{2(2s_2-s_1)}}{\pi^{2(3s_1+2s_2-1/2)}}
\ti\Gamma(s_1,s_2)^2D_{\C}(s_1,s_2)
\mho_{\C,1}(\phiu,1-s_1,1-s_2)
\end{equation*}
where the $1\times 1$ matrix $D_\C(s_1,s_2)=d_{\C,1,1}(s_1,s_2)$
is given by
\begin{equation*}
-\sin^2(\pi s_1)\sin(\pi s_2)\cos(s_1+s_2)\pi.
\end{equation*}
\end{enumerate}
\end{prop}
\begin{proof}
Let $\nabla_1,\nabla_2$, respectively be
the homogeneous linear differential operator of $y$
of degree $3,2$ with constant coefficients defined by
\begin{align*}
\nabla_1 \exp\{[y,z]_W+[\y,\z]_S\}
&	=R_1(\ti z)\exp\{[y,z]_W+[\y,\z]_S\},\\
\nabla_2 \exp\{[y,z]_W\}
&	=R_2(\ti z)\exp\{[y,z]_W\}.
\end{align*}
Then by computation we see
\begin{align*}
\nabla_1(R_1(\ti y)^{s_1}R_2(\ti y)^{s_2})
&=4s_1^2(s_1+s_2+1/2)R_1(\ti y)^{s_1-1}R_2(\ti y)^{s_2},\\
\nabla_2(R_1(\ti y)^{s_1}R_2(\ti y)^{s_2})
&=16s_2(s_1+s_2+1/2)R_1(\ti y)^{s_1}R_2(\ti y)^{s_2-1}.
\end{align*}
Hence, by the repeated use of integration by parts
we could see that
$\ti\Gamma(s_1,s_2)^{-1}\mho_{v,j_v}(\phiu_v,s_1,s_2)$
is an entire function. This proves (1).
The functional equation for the real place was
accomplished by Shintani.
The formula (2) could be found from \cite{shintanib}.
For the complex case, we choose
$\Psi(\ti y)
=\exp\{-2\pi(|y_1|_\C+2^{-1}|y_2|_\C+|y_3|_\C+|\y_1|_\C+|\y_2|_\C)\}$
as the test function.
Then as in \cite[Thorem 6.3.1]{igusa-book} the local zeta function
for $\Psi$ is computable and we could see
$\mho_{\C,1}(\Psi,s_1,s_2)
=4\pi^{-1/2}(2\pi^3)^{1-s_1}(\pi^2/4)^{1-s_2}\ti\Gamma(s_1,s_2)$.
Now we could find $c_{\C,1,1}(s_1,s_2)$ as above since $\widehat\Psi=2\Psi$.
\end{proof}
We define the product objects for infinite places as follows.
\begin{defn}
\begin{enumerate}[{\rm (1)}]
\item
We put $\type_\infty=\prod_{v\in\gM_\infty}\type_v$,
which we regard as the index set of 
orbits $H_\infty\backslash U^\sst_\infty$.
Elements of $\type_\infty$ are denoted by $j=(j_v)_{v\in\gM_\infty}$.
We put $U_j=\prod_{v\in\gM_\infty}U_{j_v}\subset U^\sst_\infty$,
which is an $H_\infty$-orbit corresponding to $j=(j_v)\in\type_\infty$.
\item
For $j\in\type_\infty$ and
$\phiu_\infty\in\sS(U_{\infty})$ we put
\begin{equation*}
\mho_{\infty,j}(\phiu_\infty,s_1,s_2)
=	\int_{U_{j}}
		|R_1(\ti y_\infty)|_\infty^{s_1-1}
		|R_2(\ti y_\infty)|_\infty^{s_2-1}
		\phiu_\infty(\ti y_\infty)d\ti y_\infty.
\end{equation*}
For $j=(j_v),l=(l_v)\in\type_\infty$, put
$d_{\infty,j,l}(s_1,s_2)=\prod_{v\in\gM_\infty}d_{k_v,j_v,l_v}(s_1,s_2)$.
If there is no confusion we drop the subscript $\infty$ and write
$\mho_j$, $d_{j,l}$ instead of $\mho_{\infty,j}$, $d_{\infty,j,l}$,
respectively.
\item
The Fourier transform $\widehat{\phiu_\infty}$
of $\phiu_\infty\in\sS(U_\infty)$ is defined similarly to
Definition \ref{defn:fourierUlocal}.
\end{enumerate}
\end{defn}
As a corollary to Proposition \ref{prop:localfneq}, we have the following:
\begin{cor}\label{cor:fneqlocal}
 For $j,l\in\type_\infty$, we have
\begin{equation*}
\mho_j(\widehat{\phiu_\infty},s_1,s_2)
=\frac{2^{n(2s_2-s_1)}}{\pi^{n(3s_1+2s_2-1/2)}}\ti\Gamma(s_1,s_2)^n
\sum_{l\in\type_\infty}
d_{j,l}(s_1,s_2)\mho_l(\phiu_\infty,1-s_1,1-s_2).
\end{equation*}
\end{cor}

\section{Density theorems}\label{sec:density}
In this section we consider the Dirichlet series
associated with the prehomogeneous vector spaces
$(G,V)$, $(B,W)$, and $(H,U)$.
By the results of Section \ref{sec:parameterization}
these turns out in Lemmata \ref{lem:zetadecompGV}, \ref{lem:zetadecompbq}
to be counting functions of cubic algebras of $\co$.
From the functional equations as well as the residue formulae
of these Dirichlet series, we derive the asymptotic formulae
on the distributions of discriminants of cubic algebras of $\co$.
Our tool to find a density theorem is a
modified version \cite[Theorem 3]{sash}
of Landau's Tauberian theorem \cite[Hauptsatz]{landaub},
using the functional equation to derive some informations on the error term.

Before starting the analysis we prepare some notation
concerning on the ideal class group $\cl(k)$ of $k$.
For a fractional ideal $\aaa$, let $i(\aaa)\in\maf(\subset \ma)$
be the corresponding idele, which is well defined up to
$\widehat\co^\times$-multiple.
That is, $i(\aaa)\in\maf$ is characterized by
the condition $\aaa=k\cap i(\aaa)\widehat\co$.
Then $|i(\aaa)|=N(\aaa)^{-1}$.
Notice that the infinite component of $i(\aaa)$ is trivial.
If there is no confusion we simply write $\aaa$ instead of $i(\aaa)$.
The set of characters of $\cl(k)$ is denoted by $\cl(k)^\ast$.
We regard $\omega\in\cl(k)^\ast$ as a character on
the idele class group $\A^1/\mk$ via the standard composition of the maps
\begin{equation*}
\A^1/\mk\rightarrow \A^1/k_\infty^1\mk\widehat\co^\times\cong 
\ma/k_\infty\mk\widehat\co^\times
\cong\maf/\mk\widehat\co^\times
	\cong \cl(k)\rightarrow \C^\times,
\end{equation*}
where we put $k_\infty^1=\A^1\cap k_\infty^\times$.
Then $\omega(\aaa)=\omega(i(\aaa))$.	
Also $\omega$ is trivial as a character on $\cl(k)$
if and only if trivial as a character on $\A^1/\mk$.
Note that this character is trivial on $\A^0$.
For an affine space $X$ defined over $\co$,
let $\sS_0(X_\A)\subset \sS(X_\A)$
be the set of elements of the form $\Phi_\infty\otimes\Phi_{\fin,0}$
where $\Phi_\infty\in\sS(X_\infty)$ is arbitrary and
$\Phi_{\fin,0}$ is the characteristic function on $X_{\widehat \co}$.
Finally, for $\omega\in\cl(k)^\ast$ let $L(s,\omega)$
be the Hecke $L$-function with respect to $\omega$;
$L(s,\omega)=\sum_{\aaa}\omega(\aaa)N(\aaa)^{-s}$
where $\aaa$ runs through the all integral ideals of $\co$.

\subsection{The Dirichlet series for the space of binary cubic forms}
We first give a remark on orbits
$G_{k_v}\backslash V^\sst_{k_v}$ for $v\in\gM_\infty$.
Let $K$ be either $\R$ or $\C$.
Then the map $\sE_2(K)\ni F\mapsto F\times K\in\sE_3(K)$ gives a bijection.
Via this map, we can construct a bijection between two sets of orbits
$H_K\backslash U^\sst_K$
and $G_K\backslash V^\sst_K$ which is useful for later purposes.
More precisely:
\begin{defn}
Let $v\in\gM_\infty$.
Then by Propositions \ref{prop:rod3} and \ref{prop:rod2},
the map $B_{k_v}\backslash W^\sst_{k_v}\to
G_{k_v}\backslash V^\sst_{k_v}$
induced by the map of Definition \ref{defn:embedding}
is bijective.
Hence we also use $\type_v$ and $\type_\infty$ as the index set
for the sets of orbits of $G_{k_v}\backslash V^\sst_{k_v}$
and $G_\infty\backslash V^\sst_\infty$.
As in Section \ref{sec:local},
the orbit corresponding to $j\in\type_v$ in $V^\sst_{k_v}$ is denoted by
$V_{k_v,j_v}$ or $V_{j_v}$, and to $j\in\type_\infty$ in $V^\sst_\infty$
is denoted by $V_j$.
\end{defn}


The Dirichlet series we are in interest are as follows.
We use the notation in Section \ref{sec:parameterization}.
\begin{defn}\label{def:zetaGV}
\begin{enumerate}[{\rm (1)}]
\item
The index set $\type_\infty$ parameterizes
the separable cubic algebra of $k_\infty$.
For $j\in\type_\infty$,
let $k_\infty(j)$ the corresponding algebra and
\begin{equation*}
\mathcal C(\mathcal O,\mathfrak a)_j
=\{R\in\mathcal C(\mathcal O,\mathfrak a)\mid
R\otimes_\co k_\infty\cong k_\infty(j)\}.
\end{equation*}
\item
We regard $\vaa$ as a lattice of $V_\infty$.
For $s\in\C$ and $\omega\in \cl(k)^\ast$
we put
\begin{align*}
\vartheta_j(\aaa;s)&=\sum_{x\in\gaa\backslash (\vaa\cap V_j)}
\frac{({}^\#(\mathrm{Stab}(\gaa;x)))^{-1}}{N(\aaa)^{2s}|P(x)|_\infty^s}
=\sum_{R\in\mathcal C(\co,\aaa)_j}
\frac{({}^\#\aut(R))^{-1}}{N(\Delta_{R/\co})^s},\\
\vartheta_j(s,\omega)&=\sum_{\aaa\in \cl(k)}\omega(\aaa)\vartheta_j(\aaa;s),
\end{align*}
where the second equality of the upper formula follows from
Proposition \ref{dfd}. Note that by this second equality,
we see that $\vartheta_j(\aaa;s)$ depends only on the ideal class of $\aaa$.
\item
For $\Phi\in\sS(V_\A)$, $s\in\C$ and $\omega\in\Omega^1$,
we define
\begin{equation*}
Z^\ast(\Phi,s,\omega)=\int_{G_\A/G_k}
|\det g|^{2s}\omega(\det g)
\sum_{x\in V^\sst_k}\Phi(gx)d_\pr g.
\end{equation*}
\item
For each $j\in\type_\infty$
we define the local zeta function at $\gM_\infty$ by
\begin{equation*}
\vcZ_{j}(\Phi_\infty,s)
=	\int_{G_{\infty}}
		|P(g_\infty x)|_\infty^{s}
		\Phi_\infty(g_\infty x)dg_\infty.
\end{equation*}
where $x$ is an arbitrary element of $V_j$.
\end{enumerate}
\end{defn}
\begin{lem}\label{lem:zetadecompGV}
For $\Phi\in\sS_0(V_\A)$, $s\in\C$ and $\omega\in\cl(k)^\ast$,
\begin{equation*}
Z^\ast(\Phi,s,\omega)=
	\sum_{j\in\mathcal T_\infty}
	\vcZ_{j}(\Phi_\infty,s)\vartheta_j(s,\omega),
\end{equation*}
\end{lem}
\begin{proof}
Let $\cK_\fin(2)=G_{\widehat\co}$.
It it is known that the double coset space
$G_\infty\cK_\fin(2)\backslash G_\A/G_k$ is represented by $\cl(k)$.
More precisely, we have
$G_\A=\coprod_{\aaa\in \cl(k)}G_\infty\cK_\fin(2)\cdot \diag(1,\aaa)\cdot G_k$.
According to this decomposition, we define the partial zeta integral by
\begin{equation*}
Z^\ast_\aaa(\Phi,s,\omega)
=	\int_{G_\infty\cK_\fin(2)\cdot \diag(1,\aaa)\cdot G_k/G_k}
		|\det g|^{2s}\omega(\det g)
\sum_{x\in V_k^\sst}\Phi(gx)d_\pr g.
\end{equation*}
We put $\cK_\fin(2)_\aaa=\diag(1,\aaa)^{-1}\cdot \cK_\fin(2)\cdot\diag(1,\aaa)$
and $\Phi_\aaa(x)=\Phi(\diag(1,\aaa)x)$.
Then since
$|\det(\cK_\fin(2)_\aaa G_k)|=\omega(\det(G_\infty\cK_\fin(2)_\aaa G_k))=1$,
$\Phi_\aaa$ is $\cK_\fin(2)_\aaa$-invariant
and $|\aaa|=N(\aaa)^{-1}$, we have
\begin{equation*}
Z^\ast_\aaa(\Phi,s,\omega)=
\frac{\omega(\aaa)}{N(\aaa)^{2s}}
\int_{G_\infty\cK_\fin(2)_\aaa/G_k\cap G_\infty\cK_\fin(2)_\aaa}
	|\det g_\infty|^{2s}_\infty\sum_{x\in V_k^\sst}\Phi_\aaa(g_\infty x)
	dg_\infty dg_\fin.
\end{equation*}
We could easily see that 
$V_k\cap \diag(1,\aaa)^{-1}V_{\widehat\co}=\vaa$
as a subset of $V_\infty$ and 
$G_k\cap \cK_\fin(2)_\aaa=\gaa$ as a subset of $G_\infty$.
Hence 
\begin{equation*}
Z^\ast_\aaa(\Phi,s,\omega)=
\frac{\omega(\aaa)}{N(\aaa)^{2s}}
\int_{G_\infty/\gaa}
	|\det g_\infty|^{2s}_\infty
	\sum_{x\in\vaa\cap V_k^\sst}\Phi_\infty(g_\infty x)
	dg_\infty
\times
\int_{\cK_\fin(2)_\aaa}dg_\fin.
\end{equation*}
Since $G_{\Af}$ is unimodular,
$\int_{\cK_\fin(2)_\aaa}dg_\fin=\int_{\cK_\fin(2)}dg_\fin=1$.
Now by the usual modification we have
\begin{equation*}
Z^\ast_\aaa(\Phi,s,\omega)=
\sum_{j\in\mathcal T_\infty}
\vcZ_j(\Phi_\infty,s)\omega(\aaa)\vartheta_j(\aaa;s).
\end{equation*}
By summing all $\aaa\in \cl(k)$ up, we have the desired formula.
\end{proof}
Hence, the Dirichlet series $\vartheta_j(s,\omega)$ is exactly
an example treated in \cite{dawra}, and
the analytic continuation, functional equation,
and residue formulae are described.
We recall the results here.
\begin{defn}\label{defn:a_kb_k}
For $j=(j_v)\in\type_\infty$, 
we define $r(j)={}^\#\{v\in\gM_\R\mid j_v=1\}$. We put
\begin{equation*}
\gA_k=\frac{\gC_k\zeta_k(2)}{2^{r_1+r_2+1}\Delta_k^{1/2}},\quad
\gB_k=
\frac{3^{r_1+r_2/2}\gC_k\zeta_k(1/3)}{5\cdot2^{r_1+r_2}\Delta_k}
\left(\frac{\Gamma(1/3)}{2\pi}\right)^{3n}.
\end{equation*}
\end{defn}
\begin{defn}\label{defn:dualGV}
Let $\hat V$ be the submodule of $V$ defined by
$\hat V=\{(x_0,3x_1,3x_2,x_3)\mid x_0,\dots,x_3\in\aff\}$.
Then $G$ acts on $\hat V$ also.
For $\aaa\in\cl(k)$ let
$\hat \vaa=V(k)\cap \diag(\mathfrak d,\mathfrak d\aaa)^{-1}\hat V(\widehat\co)$ (recall that $\mathfrak d$ is a differential idele of $k$) and put
\[
\hat\vartheta_j(s,\omega)=\sum_{\aaa\in\cl(k)}\omega(\aaa)
\sum_{x\in\gaa\backslash (\hat\vaa\cap V_j)}
{}^\#(\mathrm{Stab}(\gaa;x))^{-1}N(\aaa)^{-2s}|P(x)|_\infty^{-s}.
\]
For $j,l\in\type_\infty$,
we define $c_{j,l}(s)=\prod_{v\in\gM}c_{k_v,j_v,l_v}(s)$
where
\begin{equation*}
(c_{\R,i,j}(s))=
	\frac12\twtw{\sin2\pi s}{3\sin\pi s}{\sin\pi s}{\sin2\pi s},
\quad
c_{\C,1,1}(s)=
	\sin^2\pi s
	\sin(\pi s-\frac\pi6)
	\sin(\pi s+\frac\pi6).
\end{equation*}
\end{defn}
\begin{thm}[Datskovsky-Wright \cite{dawra}]
\label{prop:DirichletGV}
The Dirichlet series $\vartheta_j(s,\omega)$ are continued holomorphically
to the whole complex plane except for possible simple poles at $s=1,5/6$
with the residue $\delta(\omega)\gA_k(1+3^{-r(j)-r_2})$,
$\delta(\omega^3)(5/6)\gB_k3^{-r(j)/2}$,
%
respectively. They satisfy the functional equation
\begin{align*}
\vartheta_j(1-s,\omega)
=\frac{3^{n(6s-2)}}{\pi^{4ns}}\Gamma(s)^{2n}\Gamma(s-1/6)^n\Gamma(s+1/6)^n
\sum_{l\in\type_\infty}c_{l,j}(s)\hat\vartheta_l(s,\omega^{-1}).
\end{align*}
\end{thm}
\subsection{Distributions of discriminants of cubic algebras}
We are now ready to prove a density theorem from
the zeta function for $(G,V)$.
Let $\mathfrak a\in \cl(k)$.
\begin{defn}
Let $h_k^{(3)}$ be the order of the subgroup of
3-torsion elements of $\cl(k)$.
For $\mathfrak a\in \cl(k)$, we put $\tau(\mathfrak a)=1$
if $\mathfrak a$ is 3-divisible (that is, there exists
$\mathfrak b\in \cl(k)$ such that $\mathfrak a=\mathfrak b^3$),
and $\tau(\mathfrak a)=0$ otherwise.
\end{defn}
\begin{thm}\label{thm:densityGV}
\begin{enumerate}[{\rm (1)}]
\item
Let $k$ be a quadratic field. For any $\varepsilon>0$,
\begin{equation*}
\sum_{\underset{\mathcal N(\Delta_{R\slash \mathcal O})\leq X}
	{R\in \mathcal C(\mathcal O,\mathfrak a)_j}}
\frac{1}{{}^{\#}(\aut(R))}=(1+\frac1{3^{r(j)+r_2}})\frac{\mathfrak A_k}{h_k}X
+\tau(\mathfrak a)\frac{\mathfrak B_kh_k^{(3)}}{3^{r(j)/2}h_k}X^{5/6}
+O(X^{7/9+\varepsilon})\quad (X\to\infty).
\end{equation*}
\item
Let $k$ be a number field with $n=[k:\Q]\geq3$.
For any $\varepsilon>0$,
\begin{equation*}
\sum_{\underset{\mathcal N(\Delta_{R\slash \mathcal O})\leq X}
	{R\in \mathcal C(\mathcal O,\mathfrak a)_j}}
\frac{1}{{}^{\#}(\aut(R))}=(1+\frac1{3^{r(j)+r_2}})\frac{\mathfrak A_k}{h_k}X
+O(X^{(4n-1)/(4n+1)+\varepsilon})\quad (X\to\infty).
\end{equation*}
\end{enumerate}
\end{thm}
\begin{proof}
By the orthogonality of the characters,
$\vartheta_j(\aaa;s)=\sum_{\omega\in\cl(k)^\ast}
\omega(\aaa)^{-1}\vartheta(s,\omega)$.
Hence by Theorem \ref{prop:DirichletGV},
we know the analytic properties of $\vartheta_j(\aaa;s)$.
For example, the residues at $s=1,5/6$ are given by
$h_k^{-1}(1+3^{-r(j)-r_2})\mathfrak A_k$,
$\tau(\mathfrak a)(5/6)\mathfrak B_k3^{-r(j)/2}\cdot(h_k^{(3)}/h_k)$,
respectively.
Combined with \cite[Theorem 3]{sash} we have
the above formulae.
(Note that in the notation of \cite[Theorem 3]{sash},
$\nu=4n$, $\sum_{i=1}^\nu\alpha_i=0$,
$\sum_{i=1}^\nu\beta_i=4n$, $\delta=1$
and $\mu_1=\mu_2=1+\epsilon$
where $\epsilon>0$ is arbitrary.)
\end{proof}

\subsection{The Dirichlet series for the spaces of binary quadratic forms}
In this subsection we study Dirichlet series
associated with the spaces $(B,W)$ and $(H,U)$.
\begin{defn}\label{def:zetaBW}
\begin{enumerate}[{\rm (1)}]
\item
For $j\in\type_\infty$ we define $W_j\subset W^\sst_\infty$
similarly to $U_j$ and put
\begin{align*}
\xi_j(\aaa,\ccc;s_1,s_2)
&=\sum_{y\in\bac\backslash(\wac\cap W_j)}
\frac{({}^\#(\mathrm{Stab}(\bac;y)))^{-1}}{N(\ccc)^{-s_1}N(\aaa\ccc)^{2s_2}
|Q_1(y)|_\infty^{s_1}|Q_2(y)|_\infty^{s_2}},\\
\xi_j(s_1,s_2,\omega_1,\omega_2)
&=\sum_{\aaa,\ccc\in\cl(k)}
	\omega_1(\ccc)^{-1}\omega_2(\aaa\ccc)
		\xi_j(\aaa,\ccc;s_1,s_2),\\
\Xi_j(s_1,s_2,\omega_1,\omega_2)
&=\xi_j(s_1,s_2,\omega_1,\omega_2)L(2s_1,\omega_1^2).
\end{align*}
\item
We define
$Y^\ast(\phiw,s_1,s_2,\omega_1,\omega_2)
=\Delta_k^{1/2}\gC_k^2\cdot
Y(\phiw,s_1,s_2,\omega_1,\omega_2)$.
That is, we define $Y^\ast$ as the integral 
exactly the same way as in Definition \ref{zetaintBW}
except for replacing the measure $db$ by $d_\pr b$.
Similarly, we put
$X^\ast(\phiu,s_1,s_2,\omega_1,\omega_2)
=\Delta_k^{1/2}\gC_k^3\cdot
X(\phiu,s_1,s_2,\omega_1,\omega_2)$
and
$\Sigma^\ast(\cR_1\phis,s,\omega)
=\gC_k\cdot
\Sigma(\cR_1\phis,s,\omega)$.
\item
For $j\in\type_\infty$ we choose $y\in W_j$ and $\ti y\in U_j$ arbitrary.
We define 
\begin{align*}
\vcY_{\infty,j}(\phiw_\infty,s_1,s_2)
&=	\int_{B_{\infty}}
		|Q_1(b_\infty y)|_\infty^{s_1}
		|Q_2(b_\infty y)|_\infty^{s_2}
		\phiw_\infty(b_\infty y)db_\infty,\\
\vcX_{\infty,j}(\phiu_\infty,s_1,s_2)
&=	\int_{H_{\infty}}
		|R_1(h_\infty \ti y)|_\infty^{s_1}
		|R_2(h_\infty \ti y)|_\infty^{s_2}
		\phiu_\infty(h_\infty y)dh_\infty,\\
\Sigma_\infty(\cR_1\phis_\infty,s)
&=	\int_{\R^\times}|t_\infty|_\infty^s
		\phis_\infty(t_\infty,0)\, \md t_\infty.
\end{align*}
%
\end{enumerate}
\end{defn}
Recall that we put $\cK=\A^0\times\cK(2)$ which is a maximal compact
subgroup of $H_\A$.
Let $\cK_\infty$ be the infinite component of this group.
\begin{lem}\label{lem:zetadecompbq}
Assume $\phiu\in\sS(U_\A)$ is $\cK$-invariant.
We have
\begin{align*}
Y^\ast(\phiw,s_1,s_2,\omega_1,\omega_2)
&=	\sum_{j\in\mathcal T_\infty}\vcY_{\infty,j}(\phiw_\infty,s_1,s_2)
	\xi_j(s_1,s_2,\omega_1,\omega_2),\\
X^\ast(\phiu,s_1,s_2,\omega_1,\omega_2)
&=	\sum_{j\in\mathcal T_\infty}\vcX_{\infty,j}(\phiu_\infty,s_1,s_2)
	\Xi_j(s_1,s_2,\omega_1,\omega_2).
\end{align*}
\end{lem}
\begin{proof}
We put $\cK_B=B_{\widehat \co}$.
Then we have
$B_\A=\coprod_{\aaa,\ccc\in \cl(k)}B_\infty \cK_B
\cdot a(\ccc^{-1},\aaa\ccc^2)\cdot B_k$.
We could see that as a subset of $W_\infty$ or $B_\infty$,
\begin{equation*}
W_k\cap a(\ccc^{-1},\aaa\ccc^2)^{-1}W_{\widehat\co}=\wac,\quad
B_k\cap a(\ccc^{-1},\aaa\ccc^2)^{-1}\cdot \cK_B\cdot a(\ccc^{-1},\aaa\ccc^2)
=\bac.
\end{equation*}
Hence by the similar modification to
Proposition \ref{lem:zetadecompGV} we obtain the first formula.
We consider the second one.
It is well known that 
$\Sigma^\ast(\Pi',s,\omega)
=\Sigma_\infty(\Pi'_\infty,s)L(s,\omega)$
for $\Pi'\in\sS_0(\A)$.
On the other hand, since $\phiu$ is $\cK_\infty$-invariant,
we can prove
\begin{equation*}
\vcX_{\infty,j}(\phiu_\infty,s_1,s_2)
=\vcY_{\infty,j}(\phiw_\infty,s_1,s_2)\Sigma_\infty(\cR_1\phis_\infty,2s_1)
\end{equation*}
by the changes of variables similar to the proof of Lemma \ref{UW}.
Hence by Lemma \ref{UW} and the first formula we have the second one.
\end{proof}
In order to determine the analytic properties of these
Dirichlet series we recall the following well known result.
\begin{lem}\label{testfunction}
Let $j\in\type_\infty$ and $r_1,r_2\in\C$.
There exists a $\cK_\infty$-invariant function
$\phiu_\infty\in\sS(U_\infty)$ such that the support of
$\phiu_\infty$ is contained in $U_j$,
$\vcX_{\infty,j}(\phiu_\infty,s_1,s_2)$ is an
entire function and $\vcX_{\infty,j}(\phiu_\infty,r_1,r_2)\neq0$.
\end{lem}
\begin{defn}\label{defn:dualHU}
We put
$\hat U=\{(y_0,2y_1,y_2,\y_1,\y_2)\mid y_0,y_1,y_2,\y_1,\y_2\in\aff\}$,
which is a $B$-invariant submodule of $U$.
We define the Dirichlet series
$\hat\Xi_j(s_1,s_2,\omega_1,\omega_2)$
for $(B, \hat U)$ similarly to Definition \ref{def:zetaBW}.
(Also see Definition \ref{defn:dualGV}.)
\end{defn}
\begin{thm}\label{thm:fneqDHU}
The Dirichlet series
$\Xi_j(s_1,s_2,\omega_1,\omega_2)$
becomes an entire function
after multiplied by $(s_1-1)^2(s_2-1)(s_1+s_2-3/2)$.
It has the functional equation
\begin{equation*}
\Xi_j(1-s_1,1-s_2,\omega_1,\omega_2)
=\frac{2^{n(2s_2-s_1)}}{\pi^{n(3s_1+2s_2-1/2)}}\ti\Gamma(s_1,s_2)^n
\sum_{l\in\type_\infty}d_{l,j}(s_1,s_2)
\hat\Xi_j(s_1,s_2,\omega_1^{-1},\omega_2^{-1}).
\end{equation*}
\end{thm}
\begin{proof}
Take a $\cK$-invariant function $\phiu\in\sS_0(U_\A)$
such that $\phiu_\infty$ is as in Lemma \ref{testfunction}.
Then since the distributions
$\Lambda(\Psi,s_1,\omega_1)$,
$\Sigma(\cR\widehat\Psi,s_1,\omega_1)$
in Proposition \ref{pp1},
$\phis(0)$ in Proposition \ref{pp2}
and
$\Sigma_{(0)}(\cR_1\widehat{\phis},0)
\Pi_1(\phiw,s_1+2s_2,\omega_1\omega_2),
\Pi_2(\phiw\otimes\widehat{\phis},s_1+2s_2,\omega_1\omega_2)$
in Proposition \ref{pp3}
have singular supports,
$X^\ast(\phiu,s_1,s_2,\omega_1,\omega_2)$ becomes entire
after multiplied by $(s_1-1)^2(s_2-1)(s_1+s_2-3/2)$.
Also by Lemma \ref{lem:zetadecompbq},
\begin{equation*}
X^\ast(\phiu,s_1,s_2,\omega_1,\omega_2)
=\vcX_{\infty,j}(\phiu_\infty,s_1,s_2)\Xi_j(s_1,s_2,\omega_1,\omega_2).
\end{equation*}
This proves the first statement.
On the other hand since $\widehat\phiu$ is also $\cK$-invariant
and $\widehat{\phiu_\fin}$ is the characteristic function
on $\hat U_{\widehat\co}$,
the same argument of Lemma \ref{lem:zetadecompbq} shows
\begin{equation*}
X^\ast(\widehat\phiu,s_1,s_2,\omega_1,\omega_2)
=	\sum_{j\in\mathcal T_\infty}\vcX_{j}(\widehat{\phiu_\infty},s_1,s_2)
	\hat\Xi_j(s_1,s_2,\omega_1,\omega_2).
\end{equation*}
Hence by the two equations above, Theorem \ref{thm:fneqHU},
Corollary \ref{cor:fneqlocal}
and Lemma \ref{lem:b_v} we have the functional equation.
\end{proof}
%
%
\subsection{Contributions of reducible algebras}
Now we consider contributions of reducible algebras.
We describe the analytic properties of the Dirichlet series
in the question for general characters of $\cl(k)$,
but apparently the author succeed in finding density theorem
with a good error estimate
only from the Dirichlet series with trivial characters.
\begin{defn}\label{def:zetardbl}
Let $\aaa,\ccc$ be fractional ideals.
We put
\begin{align*}
\eta_j(\aaa,\ccc;s)
&=\sum_{y\in\bac\backslash (\wac\cap W_j)}
\frac{{}^\#(\mathrm{Stab}(\bac;y))^{-1}}{N(\aaa)^{2s}
|Q_1(y)^2Q_2(y)|_\infty^s},\\
\eta_j(s,\omega_1,\omega_2)
&=\sum_{\aaa,\ccc\in\cl(k)}
	\omega_1(\ccc)\omega_2(\aaa)
		\eta_j(\aaa,\ccc;s),\\
H_j(s,\omega_1,\omega_2)
&=\eta_j(s,\omega_1,\omega_2)L(4s,\omega_2^2/\omega_1^2).
\end{align*}
\end{defn}
By definition we have
$\eta_j(\aaa,\ccc;s)=\xi_j(\aaa,\ccc;2s,s)$,
$\eta_j(s,\omega_1,\omega_2)=\xi_j(2s,s,\omega_2/\omega_1,\omega_2)$
and
$H_j(s,\omega_1,\omega_2)=\Xi_j(2s,s,\omega_2/\omega_1,\omega_2)$.
The following lemma immediately follows from
Theorem \ref{thm:fneqDHU}.
\begin{lem}\label{lem:fneqzetardbl}
The series $H_j(s,\omega_1,\omega_2)$ becomes entire
after multiplied by $(s-1/2)^3(s-1)$. Moreover it satisfies
the functional equation:
\begin{multline*}
H_j(1-s,\omega_1,\omega_2)
=2^n\pi^{n(7/2-8s)}\Gamma(2s-1)^{2n}\Gamma(s)^n\Gamma(3s-3/2)^n\\
\times\sum_{l\in\type_\infty}d_{l,j}(2s-1,s)
\hat\Xi_j(2s-1,s,\omega_1/\omega_2,1/\omega_2).
\end{multline*}
\end{lem}
We next consider the residue at $s=1$.
\begin{lem}\label{lem:residuezetardbl}
The residue of $\eta_j(s,\omega_1,\omega_2)$ at $s=1$
is $\delta(\omega_2)\mathfrak A_k L(2,\omega_1^{-1})/\zeta_k(2)$.
\end{lem}
\begin{proof}
Take $\phiw\in\sS_0(W_\A)$ such that $\supp(\phiw_\infty)\in W_j$.
By Lemma \ref{lem:zetadecompbq} we have
\begin{equation*}
Y^\ast(\phiw,2s,s,\omega_2/\omega_1,\omega_2)
=\vcY_{\infty,j}(\phiw_\infty,2s,s)\eta_j(s,\omega_1,\omega_2)
\end{equation*}
We consider the residue at $s=1$. By Proposition \ref{pp1} we have
\begin{align*}
\res_{s=1}Y^\ast(\phiw,2s,s,\omega_2/\omega_1,\omega_2)
&=\delta(\omega_2)2^{-1}\Delta_k^{1/2}\gC_k^2
	\Sigma(\cR\phiw,2,\omega_2/\omega_1)\\
&=\delta(\omega_2)2^{-1}\Delta_k^{1/2}\gC_k^2
	\int_{\ma\times\A^2}\omega_1^{-1}(t)|t|^2\phiw(t,u_2,u_3)
		\,\md tdu_1du_2\\
&=\delta(\omega_2)2^{-1}\Delta_k^{-1/2}\gC_kL(2,\omega_1^{-1})
	\int_{k_\infty^3}|y_{1,\infty}|_\infty
		\phiw_\infty(y_\infty)dy_\infty.
\end{align*}
On the other hand for $v\in\gM_\infty$ we could see
\begin{equation*}
\int_{B_{k_v}}f(b_vy)db_v=
2\int_{B_{k_v}y}f(z_v)|Q_1(z_v)|_v^{-1}|Q_2(z_v)|_v^{-1} dz_v
\end{equation*}
by computing the Jacobian determinant of the double cover
$B_{k_v}\ni b_v\mapsto b_vy\in B_{k_v}y$ for their coordinate systems.
This shows that
\begin{equation*}
\vcY_{\infty,j}(\phiw_\infty,2,1)
=2^{r_1+r_2}\int_{W_j}|R_1(y_{1,\infty})|_\infty\phiw_\infty(y_\infty)dy_\infty
=2^{r_1+r_2}\int_{W_\infty}|y_{1,\infty}|_\infty\phiw_\infty(y_\infty)dy_\infty
\end{equation*}
since $\supp(\phiw_\infty)\subset W_j$.
Hence we have the proposition.
\end{proof}
To find the density theorem we prepare a lemma.
\begin{lem}\label{lem:zeta^-1}
The Dirichlet series $1/\zeta_k(s)=\prod_{v\in\gMf}(1-N(\mathcal P_v)^{-s})$
is normally convergent for $\Re(s)>1$,
where we denote by $\mathcal P_v$ by the prime ideal of $\co$
corresponding to $v\in\gMf$.
\end{lem}
\begin{proof}
For two Dirichlet series
$a(s)=\sum_{n\geq1}a_n/n^s$ and $b(s)=\sum_{n\geq1}b_n/n^s$
we write $a(s)\preceq b(s)$ to mean that $|a_n|\leq b_n$ for all $n$.
Let $p$ be a prime.
For $k\in\mathbb N_{\geq1}$ we have
$1-p^{-ks}\preceq(1+p^{-s})^k\preceq(1-p^{-s})^{-k}$.
Hence $\prod_{v\mid p}(1-N(\mathcal P_v)^{-s})\preceq(1-p^{-s})^{-n}$.
Taking the product for all primes,
we have $\zeta_k(s)^{-1}\preceq\zeta_\Q(s)^n$.
Since the right hand side is normally convergent for $\Re(s)>1$,
the left hand side also.
\end{proof}
\begin{prop}\label{prop:rdbldensity}
For any $\varepsilon>0$,
\begin{equation*}
\sum_{\aaa,\ccc\in\cl(k)}
\sum_{\underset{N(\aaa)|Q_1(y)^2Q_2(y)|_\infty<X}
{y\in\bac\backslash(\wac\cap W_j)}}
\frac{1}{{}^\#(\mathrm{Stab}(\bac;y))}
	=\mathfrak A_kX+O(X^{(5n-1)/(5n+1)+\varepsilon})
\quad(X\to\infty).
\end{equation*}	
\end{prop}
\begin{proof}
We write $\eta_j(s)=\eta_j(s,1,1)$ and $H_j(s)=H_j(s,1,1)$,
where $1\in\cl(k)^\ast$ is the trivial character.
We fix $j\in\type_\infty$.
Let $\eta_j(s)=\sum_{n\geq1} a_n/n^s$,
so that the left hand side of the proposition
is $\sum_{n<X}a_n$. We denote this function by $A(X)$.
We also put $H_j(s)=\sum_{n\geq1} b_n/n^s$ and
$B(X)=\sum_{n<X}b_n$.
By Definition \ref{def:zetardbl} and Lemma \ref{lem:residuezetardbl}
$H_j(s)$ is holomorphic for $\Re(s)>1/2$
except for a simple pole at $s=1$
with the residue $\mathfrak A_k\zeta_k(4)$.
Using the functional equation in Lemma \ref{lem:fneqzetardbl},
by \cite[Theorem 3]{sash}
we have $B(X)=\mathfrak A_k\zeta_k(4)X+O(X^{(5n-1)/(5n+1)+\varepsilon})$.
Note that in the notation of \cite[Theorem 3]{sash},
$\nu=4n$, $\sum_{i=1}^\nu\alpha_i=-7n/2$,
$\sum_{i=1}^\nu\beta_i=8n$, $\delta=1$
and $\mu_1=\mu_2=1+\epsilon$
where $\epsilon>0$ is arbitrary.

We put $1/\zeta_k(4s)=\sum_{n\geq1} d_n/n^s$.
Then by Definition \ref{def:zetardbl}, $a_n=\sum_{mk=n}d_mb_k$.
We define $\rho(X)=B(X)-\mathfrak A_k\zeta_k(4)X$,
$\mu=(5n-1)/(5n+1)$.
Then for any $\varepsilon>0$, 
there exists $M>0$ such that $|\rho(X)|\leq MX^{\mu+\varepsilon}$.
Hence we have
\begin{align*}
A(X)
&=	\sum_{mk<X}d_mb_k
=	\sum_{m<X}d_m\sum_{k<X/m}b_k
=	\sum_{m<X}d_mB(X/m)\\
&=	\mathfrak A_k\zeta_k(4)X\sum_{m<X}\frac{d_m}{m}
+	\sum_{m<X}d_m\rho(X/m)\\
&=	\mathfrak A_k\zeta_k(4)X
		\left(\frac{1}{\zeta_k(4)}-\sum_{m\geq X}\frac{d_m}{m}\right)
+	\sum_{m<X}d_m\rho(X/m).
\end{align*}
Since
\begin{align*}
\sum_{m\geq X}\frac{|d_m|}{m}
&\leq	\sum_{m\geq X}\frac{|d_m|}{m^{1/4+\varepsilon}}X^{-3/4+\varepsilon}
=	O(X^{-3/4+\varepsilon}),\\
\sum_{m<X}|d_m||\rho(X/m)|
&\leq	\sum_{m<X}|d_m|\cdot M\frac{X^{\mu+\varepsilon}}{m^{\mu+\varepsilon}}
=	MX^{\mu+\varepsilon}\sum_{m<X}\frac{|d_m|}{m^{\mu+\varepsilon}}
=	O(X^{\mu+\varepsilon}),
&\end{align*}
we have the proposition.
\end{proof}
\subsection{The distributions of irreducible algebras}
We are now ready to prove a main theorem of this paper.
\begin{defn}
\begin{enumerate}[{\rm (1)}]
\item
We denote by $L_\infty$ a separable cubic $k_\infty$-algebra.
Let $r(L_\infty)$ be the integer
such that $L_\infty\cong \R^{r_1+2r(L_\infty)}\times\C^{3r_2+r_1-r(L_\infty)}$
as $\R$-algebras.
Notice that $0\leq i(L_\infty)\leq r_1$.
\item
Let $h(n,L_\infty)$ be the numbers of
the isomorphism classes of pairs $(R,F)$ satisfying the following conditions:
(a) $F/k$ is a cubic (field) extension in $\overline\Q$ such that
$F\otimes_{k}k_\infty\cong L_\infty$ as $k_\infty$-algebras,
(b) $R$ is an order of $F$ containing $\co$,
(c) $N(\Delta_{R/\co})=n$.
\end{enumerate}
\end{defn}
Note that $N(\Delta_{R/\co})$ is the ideal norm of
the relative discriminant of $R/\co$.
The following is a main result of this paper
\footnote{We define constants $\mathfrak A_k,\mathfrak B_k$
in Definition \ref{defn:a_kb_k}.}.
\begin{thm}\label{thm:field}
\begin{enumerate}[{\rm (1)}]
\item
Let $k$ be a quadratic field. For any $\varepsilon>0$,
\begin{equation*}
\sum_{n<X}h(n,L_\infty)
=3^{-r(L_\infty)-r_2}\mathfrak A_kX
+3^{-r(L_\infty)/2}\mathfrak B_kX^{5/6}
+O(X^{9/11+\varepsilon})\quad(X\to\infty).
\end{equation*}
\item
Let $k$ be a number field with $n=[k:\Q]\geq3$.
For any $\varepsilon>0$,
\begin{equation*}
\sum_{n<X}h(n,L_\infty)
=3^{-r(L_\infty)-r_2}\mathfrak A_kX
+O(X^{(5n-1)/(5n+1)+\varepsilon})\quad(X\to\infty).
\end{equation*}
\end{enumerate}
\end{thm}
\begin{proof}
Let $\mathcal C(\co)_j=\coprod_{\aaa\in\cl(k)}\mathcal C(\co,\aaa)_j$.
Since $\sum_{\aaa\in\cl(k)}\tau(\aaa)=h_k/h_k^{(3)}$,
by adding all $\aaa\in\cl(k)$ up
for the formula of Theorem \ref{thm:densityGV} we have
\begin{equation*}
\sum_{\underset{\mathcal N(\Delta_{R\slash \mathcal O})\leq X}
	{R\in \mathcal C(\mathcal O)_j}}
\frac{1}{{}^{\#}(\aut(R))}=(1+\frac1{3^{r(j)+r_2}}){\mathfrak A_k}X
+{\mathfrak B_k3^{-r(j)/2}}X^{5/6}
+O(X^{(4n-1)/(4n+1)+\varepsilon}).
\end{equation*}
We will compare the formula of Proposition \ref{prop:rdbldensity}
and the formula above.
Since $\aut_\co(R)$ is a subgroup of $\aut_k(R\otimes_\co k)$,
the order of $\aut_\co(R)$ is either $1,2,3$ or $6$.
Especially if $R\otimes_\co k$ is a non-cyclic cubic extension,
then $\aut_\co(R)$ is trivial.
We consider the algebras $R$ such that
$\aut_k(R\otimes_\co k)$ is isomorphic
to either $\Z/3\Z$ or $S_3$, the permutation group of degree $3$.
We denote by ${\rm CG}_k$ the set of cyclic cubic extensions of $k$.
Let $p(n)$ be the number of orders $R$ of $k\times k\times k$
with $N(\Delta_{R/\co})=n$, and $q(n)$ the the numbers of
the pairs $(R,F)$ such that $F\in{\rm CG}_k$,
$R$ is an order of $F$, and $N(\Delta_{R/\co})=n$.
We claim that the contributions of these algebras
can be ignored in the limit. To see this, we
notice that for a positive sequence $\{a_n\}$ and
a positive constant $\rho$,
the series $\sum_{n\geq1}a_n/n^s$ converges for $\Re(s)>\rho$
if and only if
$\sum_{n<X}a_n=O(X^{\rho+\epsilon})$
for any $\epsilon>0$.
By \cite[Theorem 6.1]{dawra} we have
\begin{align*}
\sum_{n\geq1}\frac{p(n)}{n^s}
&=\zeta_k(2s)^3\zeta_k(6s-1)/\zeta_k(4s)^2,\\
\sum_{n\geq1}\frac{q(n)}{n^s}
&=\sum_{F\in{\rm CG}_k}\frac{1}{N(\Delta_{F/k})^s}
\zeta_k(4s)\zeta_k(6s-1)\zeta_F(2s)/\zeta_F(4s).
\end{align*}
Let $\epsilon$ be an arbitrary positive number.
From the first equality, we have $\sum_{n<X}p(n)=O(X^{1/3+\epsilon})$.
Also by \cite[Theorem I.2]{wrightb} we know that
$\sum_{F\in{\rm CG}_k}N(\Delta_{F/k})^{-s}$
converges for $\Re(s)>1/2$.
Since the Dirichlet series $\zeta_F(2s)/\zeta_F(4s)$
is uniformly bounded by $\zeta_k(2s)^3/\zeta_k(4s)^3$,
the Dirichlet series in the right hand side of the second formula
converges for $\Re(s)>1/2$,
which asserts $\sum_{n<X}q(n)=O(X^{1/2+\epsilon})$.
Now the theorem follows from Propositions
\ref{prop:reducibleparameterization} and \ref{prop:rdbldensity}.
\end{proof}

\begin{rem}
Let us consider the case $k=\Q$.
We see that $\mathfrak A_\Q=\pi^2/24$ and $\mathfrak B_\Q=r/10$
where we put $r=(2\pi)^{1/3}\zeta(2/3)\Gamma(1/3)\Gamma(2/3)^{-1}$.
Let $h(n)$ be the numbers of orders of
isomorphism classes of cubic fields with discriminant $n$.
Then the proof of Theorem \ref{thm:field} also shows
\begin{align*}
\sum_{0<n<X}h(n)&=(\pi^2/72)X+
	(\sqrt3r/30)X^{5/6}+O(X^{2/3+\epsilon})\quad(X\to\infty),\\
\sum_{0<n<X}h(-n)&=(\pi^2/24)X+
	(r/10)X^{5/6}+O(X^{2/3+\epsilon})\quad(X\to\infty).
\end{align*}
This is what Shintani established in \cite[Theorem 4]{shintanib}
and hence Theorem \ref{thm:field} is a generalization of this to
an arbitrary number field.
Note that Shintani used $\spl(2)_\Z$ instead of $\gl(2)_\Z$
and hence his result is twice to ours.
\end{rem}

\section{On zeta integrals for the space of binary cubic forms}
\label{sec:cubic}
In this section we consider some zeta integrals
concerning on the space of binary cubic forms $(G,V)$.
We give meromorphic continuations and
describe some of the residues of these integrals
in Proposition \ref{prop:zetaintegral}.

Recall that we associate $k(x)\in\sE_3(k)$
for each $x\in V_k^\sst$ in Definition \ref{def:rod3}.
Let
\begin{align*}
V_{(1)}&=\{x\in V_k^\sst\mid k(x)\cong k\times k\times k\},\\
V_{(2)}&=\{x\in V_k^\sst\mid k(x)\cong k\times E
			\text{ where $E/k$ is a quadratic field extension}\},\\
V_{(3)}&=\{x\in V_k^\sst\mid
			\text{$k(x)/k$ is a cubic field extension}\}.
\end{align*}
Then $V_k^\sst=V_{(1)}\amalg V_{(2)}\amalg V_{(3)}$
and each of them is a $G_k$-invariant subset.
\begin{defn}
For $\Phi\in\sS(V_\A)$, $s\in\C$ and $\omega\in\Omega^1$,
the global zeta function is defined by
\begin{equation*}
Z(\Phi,s,\omega)=\int_{G_\A/G_k}|\det g|^{2s}\omega(\det g)
\sum_{x\in V^\sst_k}\Phi(gx)dg.
\end{equation*}
For $i=1,2,3$, we also define the {\em zeta integrals} for $V_{(i)}$ by
\begin{equation*}
Z_{(i)}(\Phi,s,\omega)=\int_{G_\A/G_k}|\det g|^{2s}\omega(\det g)
\sum_{x\in V_{(i)}}\Phi(gx)dg.
\end{equation*}
\end{defn}
By definition
$Z(\Phi,s,\omega)=\sum_{1\leq i\leq3}Z_{(i)}(\Phi,s,\omega)$.
We give some analytic properties of integrals $Z_{(i)}(\Phi,s,\omega)$.
Note that $V_{(1)}$ is a single $G_k$-orbit and
the results of $Z_{(1)}(\Phi,\omega,s)$ in this section
is included in Datskovsky-Wright's analysis \cite{dawra}
of the orbital zeta functions.
Before starting the analysis,
without loss of generality we assume the following
as in Section \ref{sec:global}.

\begin{asmp}
The Schwartz-Bruhat function satisfies $\cM_\omega\Phi=\Phi$,
where the operator $\cM_\omega$ is defined by
$\cM_\omega\Phi(x)=\int_{\cK(2)}\omega(\det\kappa)\Phi(\kappa x)d\kappa$.
\end{asmp}
We will express $Z_{(i)}(\Phi,s,\omega)$ $(i=1,2)$
by means of following zeta integrals of $(B,W)$.
\begin{defn}
Let $W_{(1)}=\{y\in W_k^\sst\mid k(y)\cong k\times k\}$
and $W_{(2)}=W_k^\sst\setminus W_{(1)}$.
For $i=1,2$, we define the zeta integrals of $(B,W)$ by
\begin{equation*}
Y_{(i)}(\phiw,s_1,s_2,\omega_1,\omega_2)
 =	\int_{B_\A/B_k}
	|t|^{s_1}|tp|^{2s_2}
		\omega_1(t)\omega_2(tp)
			\sum_{y\in W_{(i)}}\phiw(by)db.
\end{equation*}
\end{defn}
\begin{lem}\label{lem:ZY}
Let us define $\cR_W\Phi\in\sS(W_\A)$ by
$\cR_W\Phi(y)=\Phi(y^\ast)$.
We have
\begin{align*}
Z_{(2)}(\Phi,s,\omega)
&	=Y_{(2)}(\cR_{W}\Phi,2s,s,\omega,\omega),\\
Z_{(1)}(\Phi,s,\omega)
&	=3^{-1}Y_{(1)}(\cR_{W}\Phi,2s,s,\omega,\omega).
\end{align*}
\end{lem}
\begin{proof}
We can easily check that $V_{(2)}=G_k\times_{\br(2)_k}(W_{(2)})^\ast$.
Hence
\begin{align*}
Z_{(2)}(\Phi,s,\omega)
&=	\int_{G_\A/\br(2)_k}|\det g|^{2s}\omega(\det g)
	\sum_{y\in W_{(2)}}\Phi(gy^\ast)dg\\
&=	\int_{\br(2)_\A/\br(2)_k}|\det b_2|^{2s}\omega(\det b_2)
	\sum_{y\in W_{(2)}}\Phi(b_2y^\ast)db_2\\
&=	\int_{B_\A/B_k}|t^2p|^{2s}\omega(t^2p)
	\sum_{y\in W_{(2)}}\Phi(b^\ast y^\ast)db
\end{align*}
which is equal to $Y_{(2)}(\cR_{W}\Phi,2s,s,\omega,\omega)$
since $\Phi(b^\ast y^\ast)=\cR_W\Phi(by)$.
This proves the first equality.
Let $w=(1,1,0)\in W_{(1)}$.
Then $w^\ast=(0,1,1,0)\in V_{(1)}$.
Since $V_{(1)}$ is a single $G_k$-orbit and
${\rm Stab}(G_k;w^\ast)$ is of order $6$, we have
\begin{equation*}
Z_{(1)}(\Phi,s,\omega)
=	6^{-1}\int_{B(2)_\A/B(2)_k}|\det b_2|^{2s}\omega(\det b_2)
	\Phi(b_2w^\ast)db_2.
\end{equation*}
On the other hand,
$W_{(1)}$ is also a single $B_k$-orbit and
${\rm Stab}(B_k;w)$ is of order $2$.
Hence by a similar modification we obtain the second formula.
\end{proof}
\begin{defn}
For $\Phi\in\sS(V_\A)$, we define
$\cR_a\Phi,\cR_b\Phi\in\sS(\A)$ by
\begin{equation*}
\cR_a\Phi(x)=\int_{\A^2}\Phi(0,x,y_1,y_2)dy_1dy_2,\quad
\cR_b\Phi(x)=\int_{\A^3}\Phi(x,y_1,y_2,y_3)dy_1dy_2dy_3.
\end{equation*}
\end{defn}
\begin{prop}\label{prop:zetaintegral}
\begin{enumerate}[{\rm (1)}]
\item
The zeta integral $Z_{(1)}(\Phi,s,\omega)$ is holomorphic for $\Re(s)>1/3$.
\item
The zeta integral $Z_{(2)}(\Phi,s,\omega)$ is holomorphic for $\Re(s)>1/2$
except for possible simple pole at $s=1$ with the residue
$\delta(\omega)2^{-1}\Sigma(\cR_a\Psi,2)$.
\item
The zeta integral $Z_{(3)}(\Phi,s,\omega)$ is holomorphic for $\Re(s)>1/2$
except for possible simple poles at $s=1,5/6$ with the residue
$\delta(\omega)2^{-1}\gC_k^{-1}\Delta_k\zeta_k(2)\int_{V_\A}\Phi(x)dx$,
$\delta(\omega^3)6^{-1}\Sigma(\cR_b\Phi,1/3)$,
respectively.
\end{enumerate}
Moreover each of the zeta integrals has meromorphic continuation
to the whole complex plane.
\end{prop}
\begin{proof}
(1) follows from either \cite[Theorem 6.1]{dawra} or
Lemma \ref{lem:ZY} and Proposition \ref{prop:orbitalzeta},
and (2) follows from (1), Lemma \ref{lem:ZY}, Proposition \ref{pp1}
and Theorem \ref{thm:fneqHU}.
On the other hand the analytic properties including
meromorphic continuations and residue formulae for
$Z(\Phi,s,\omega)$ are known by \cite[Theorem 6.4]{wright}.
(3) follows from this and (1), (2).
The meromorphic continuations are also proved
 one by one from (1) to (3).
\end{proof}

\section{Acknowledgments}
During this work, T. Terasoma gave the author useful suggestions.
Especially the idea of Proposition \ref{prop:reducibleparameterization}
was suggested by him.
The author express his gratitude to him.
The author is also grateful to A. Kable and D. Wright
for their communications as well as making
\cite{kawr} available to him prior to its publication,
which becomes the starting point of this work.
The author also thanks to F. Sato and K. Matsumoto for helpful communications.
Financial support was provided by
21st Century (the University of Tokyo) COE program,
of the Ministry of Education, Culture, Sports, Science and Technology.

\appendix
\section{Remaining functional equations}\label{sec:remainfneq}

Here, we collect $2$ kinds of local functional equations for $(H,U)$.
For the real case, the first equation
was proved in \cite[Lemma 1.9]{shintania}
and the second one is proved in \cite[Lemma 2.9]{fsatoa}.
The complex case immediately follows from
the proof of Proposition \ref{prop:localfneq}.
We note that the functional equations in
Proposition \ref{prop:localfneq} follows from this table also.

\begin{prop}
\begin{enumerate}[{\rm (1)}]
\item The functional equations for $\R$ are:
\begin{multline*}
\begin{pmatrix}
\mho_{\R,1}(\widehat\phiw\otimes\phis,s_1,s_2)\\
\mho_{\R,2}(\widehat\phiw\otimes\phis,s_1,s_2)\\
\end{pmatrix}
=\frac{2^{s_1+2s_2-1}}{\pi^{s_1+2s_2-1/2}}
\Gamma(s_2)\Gamma(s_1+s_2-1/2)\\
\times
\begin{pmatrix}
\sin((s_1+2s_2)\pi/2)	&\cos(s_1\pi/2)\\
\sin(s_1\pi/2)		&\cos((s_1+2s_2)\pi/2)\\
\end{pmatrix}
\begin{pmatrix}
\mho_{\R,1}(\phiw\otimes\phis,s_1,3/2-s_1-s_2)\\
\mho_{\R,2}(\phiw\otimes\phis,s_1,3/2-s_1-s_2)\\
\end{pmatrix},
\end{multline*}
\begin{multline*}
\begin{pmatrix}
\mho_{\R,1}(\phiw\otimes\widehat\phis,s_1,s_2)\\
\mho_{\R,2}(\phiw\otimes\widehat\phis,s_1,s_2)\\
\end{pmatrix}
=\frac{2}{(2\pi)^{2s_1}}
\Gamma(s_1)^2\\
\times
\begin{pmatrix}
2\cos^2(s_1\pi/2)	&0\\
0			&\sin(s_1\pi)\\
\end{pmatrix}
\begin{pmatrix}
\mho_{\R,1}(\phiw\otimes\phis,1-s_1,s_1+s_2-1/2)\\
\mho_{\R,2}(\phiw\otimes\phis,1-s_1,s_1+s_2-1/2)\\
\end{pmatrix}.
\end{multline*}

\item The functional equations for $\C$ are:
\begin{multline*}
\mho_{\C,1}(\widehat\phiw\otimes\phis,s_1,s_2)
=\frac{2^{2s_1+4s_2-2}}{\pi^{2s_1+4s_2-1}}
\Gamma(s_2)^2\Gamma(s_1+s_2-1/2)^2\\
\times\sin(s_2\pi)\cos(s_1\pi+s_2\pi)
\mho_{\C,1}(\phiw\otimes\phis,s_1,3/2-s_1-s_2),
\end{multline*}
\begin{multline*}
\mho_{\C,1}(\phiw\otimes\widehat\phis,s_1,s_2)
=\frac{4}{(2\pi)^{4s_1}}
\Gamma(s_1)^4\sin^2(s_1\pi)
\mho_{\C,1}(\phiw\otimes\phis,1-s_1,s_1+s_2-1/2).
\end{multline*}
\end{enumerate}
\end{prop}

By the same arguments of Section \ref{sec:density}
we could find $2$ more functional equations for
Dirichlet series $\Xi_j(s_1,s_2,\omega_1,\omega_2)$.

\section{Non-archimedean local theory for the spaces of binary quadratic forms}
\label{sec:non-arch}
Here we study some non-archimedean local theory.
The explicit formula for the standard test function is obtained.
An interesting corollary to this formula, we describe
the orbital zeta functions by means of Dedekind zeta functions.

\subsection{Explicit formula at finite places}
In this subsection we give the explicit formula
of the local zeta function for the standard test function
at finite places.

For a while let $v\in\gM$ arbitrary (not assuming a finite place).
To treat ramified characters also,
we slightly generalize the notion of local zeta function.
\begin{defn}\label{localzeta}
For $\ti y\in U_{k_v}^\sst$, $\phiu_v\in \sS(U_{k_v})$,
$\omega_{1v},\omega_{2v}\in\Omega_v$ and $s_1,s_2\in\C$,
we define
\begin{equation*}
\gpX_{\ti y,v}(\phiu_v,s_1,s_2,\omega_{1v},\omega_{2v})
 =\int_{H_{k_v}} 
	\omega_{1v}(t_v)
	\omega_{2v}(t_v\det g_v)
	|t_v|_v^{s_1}|t_v\det g_v|_v^{2s_2}
	\Phi_v(h_v\ti y)\, dh_v
\end{equation*}
and call it {\em the local zeta function}.
\end{defn}
If both $\omega_{1v},\omega_{2v}$ are trivial
we often drop it and write $\gpX_{\ti y,v}(\phiu_v,s_1,s_2)$.
For $j_v\in\type_v$ such that $\ti y\in U_{j_v}$,
by definition
$\vcX_{v,j_v}(\phiu_v,s_1,s_2)=
|R_1(\ti y)|_v^{s_1}|R_2(\ti y)|_v^{s_2}
\gpX_{\ti y,v}(\phiu_v,s_1,s_2)$.

For any $\ti y\in U_{k_v}^\sst$,
The stabilizer of $\ti y$ in $H_{k_v}$
consists of two elements
and the non-trivial element is of the form $(1,g)$
with $\det(g)=-1$.
This shows that $\gpX_{\ti y,v}(\phiu_v,s_1,s_2,\omega_{1v},\omega_{2v})$
is identically zero unless $\omega_{2v}(-1)=1$.
Hence we assume the following.
\begin{asmp}\label{omega2}
For any $v\in\gM$, $\omega_{2v}(-1)=1$.
\end{asmp}

The analytic continuation of the local zeta function
is known in more general settings than the prehomogeneous case.
The meromorphic properties of complex powers of polynomials were
studied by Bernstein and Gelfand \cite{bege} for infinite places
and by Denef \cite{denefa, denefb} for finite places.
The following lemma is contained in their works.
(We gave a proof for $v\in\gM_\infty$
in Proposition \ref{prop:localfneq}.)
\begin{lem}\label{localcont}
The local zeta function
$\gpX_{\ti y,v}(\phiu_v,s_1,s_2,\omega_1,\omega_2)$
has meromorphic continuation to $\C^2$.
Moreover, it is a rational function
of $q_v^{s_1},q_v^{2s_2}$ if $v\in\gMf$.
\end{lem}

It will be convenient to
attach to each orbit in $U_{k_v}^{\sst}$ 
where $v\in\gM$, an index or type which
records the arithmetic properties of $v$ and the extension of
$k_v$ corresponding to the orbit.
Recall that by Proposition \ref{prop:rod2},
the orbit space $H_{k_v}\backslash U_{k_v}^\sst$
corresponds bijectively to the set of 
isomorphism classes of separable quadratic algebra of $k_v$.
The orbit corresponding to
$k_v\times k_v$ will have the index (sp).
The orbit corresponding to the unique unramified quadratic extension of
$k_v$ will have the index (ur).
An orbit corresponding to a ramified quadratic extension of $k_v$
will have the index (rm).
Recall that the extension $\C/\R$ is regarded as ramified.

\begin{defn}\label{standardrepresentative}
For each of $H_{k_v}$-orbits in $U_{k_v}^\sst$,
we choose and fix an element $\ti z$ which satisfies the following condition.
\begin{enumerate}[(1)]
\item If the orbit is corresponding to $k_v\times k_v$,
then $\ti z=(0,1,0,1,1)$.
\item If $v\in\gM_\R$ and the orbit is corresponding to $\C$,
then $\ti z=(1/2,0,1/2,1,1)$.
\item Consider the case $v\in\gMf$ and the orbit is corresponding
to a quadratic field extension $F/k_v$.
We choose and fix $z\in W_{k_v}$ such that $z_1=1$ and
$\co_F$ is generated over $\co_v$ by roots of $z(v_1,1)$.
We put $\z=(1,0)$ and choose $\ti z=(z,\z)$ as the orbital representative.
\end{enumerate}
We call such fixed orbital representatives as
the {\em standard orbital representatives}.
\end{defn}

We note that for any standard representative $\ti z$, we have $R_1(\ti z)=1$,
and if $v\in\gMf$ then the discriminant $R_2(\ti z)$ of $z(v)$
generates the ideal $\Delta_{\ti k_v(\ti z)/k_v}$.

We now assume $v\in\gMf$.
We give an explicit formula when $\phiu_0$
is the characteristic function of $U_{\co_v}$.
In this case the integral is $0$
either $\omega_{1v}$ or $\omega_{2v}$ is ramified.
Hence we consider the case
both $\omega_{1v}$ and $\omega_{2v}$ are unramified.

\begin{prop}\label{explicitformula}
For any $v\in\gMf$, let $\phiu_{v,0}$ be
the characteristic function of $U_{\co_v}$.
For a standard representative $\ti z\in U_{k_v}^\sst$,
\begin{equation*}
\gpX_{\ti z,v}(\phiu_{v,0},s_1,s_2)
=(1-q_v^{-2s_2})^{-1}(1-q_v^{1-2s_1-2s_2})^{-1}(1-q_v^{-s_1-2s_2})^{-1}
R_{v,\ti z}(s_1,s_2)
\end{equation*}
where
\begin{equation*}
R_{v,\ti z}(s_1,s_2)
=
\begin{cases}
(1-q_v^{-s_1-2s_2})^2/(1-q_v^{-s_1})^{2}&
	\text{$\ti z$ is of type {\rm (sp)}},\\
(1-q_v^{-2s_1-4s_2})/(1-q_v^{-2s_1})&
	\text{$\ti z$ is of type {\rm (ur)}},\\
(1-q_v^{-s_1-2s_2})/(1-q_v^{-s_1})&
	\text{$\ti z$ is of type {\rm (rm)}}.\\
\end{cases}
\end{equation*}
\end{prop}
\begin{proof}
In the proof of this proposition,
we drop the subscript $v$ from various symbols
such as $\gpX_{\ti z,v}, \phiu_{0,v}, u_v, q_v$ if there is no confusion.
We first consider the case $\ti z$ is of type (sp).
We put $a=s_1+2s_2, b=s_1$.
Then, by a standard modification, we have
\begin{equation*}
\gpX_{\ti z}(\phiu_{0},s_1,s_2)
=	\int_{(\mk_v)^3\times k_v}|t_1|_v^a|t_2t_3|_v^b
		\phiu_0(0,t_1,t_1u,t_2-ut_3,t_3)
		\md t_1\md t_2\md t_3 du.
\end{equation*}
Let $\pi\in\co_v$ be a uniformizer.
By changing $t_1$ to $\pi^{-1} t_1$,
$t_3$ to $\pi^{-1} t_3$, and $u$ to $\pi u$,
we have
\begin{equation*}
\gpX_{\ti z}(\phiu_{0},s_1,s_2)
=q^{a+b-1}
	\int_{(\mk_v)^3\times k_v}|t_1|_v^a|t_2t_3|_v^b
		\phiu_0(0,t_1/\pi,t_1u,t_2-ut_3,t_3/\pi)
		\md t_1\md t_2\md t_3 du.
\end{equation*}
Hence if we let $\phiu_0'\in\sS(U_{k_v})$
be the characteristic function
of $U_{\co_v}\setminus(\co_v\times\gp_v\times\co_v\times\co_v\times\gp_v)$,
then
$\gpX_{\ti z}(\phiu_0,s_1,s_2)
	=(1-q^{1-a-b})^{-1}\gpX_{\ti z}(\phiu_0',s_1,s_2)$.
We consider the integral $\gpX_{\ti z}(\phiu_0',s_1,s_2)$.
We divide the domain of the integration into the following three subsets:
(a) $t_1\in\co_v^\times$,
(b) $t_1\in\pi^m\co_v^\times, t_3\in\co_v^\times, u\in\co_v$ for $m\geq1$,
(c) $t_1\in\pi^m\co_v^\times, t_3\in\co_v^\times, u\in(\pi^{-m}\co_v\setminus\co_v)$ for $m\geq1$.
Then the value of the integral in each domain is found to be
\begin{equation*}
(1-q^{-b})^{-2},
\quad
q^{-a}(1-q^{-a})^{-1}(1-q^{-b}),
\quad
\text{and}
\quad
q^{-a+b}(1-q^{-a})^{-1}(1-q^{-a+b}),
\end{equation*}
respectively.
Adding all these up, we have
$\gpX_{\ti z}(\phiu_0',s_1,s_2)=(1-q^{-a})(1-q^{-b})^{-2}(1-q^{-a+b})^{-1}$,
and this proves the formula for the case (sp).

We next consider the case $x$ is of type (ur) or (rm).
We put $F=\ti k_v(\ti z)$.
Let $\cT,\cN$ be the trace and the norm of the quadratic extension $F/k$.
We write $x_1(v)=\cN(v_1+\theta v_2)$ where $\theta\in\co_F$.
Then by the definition of the standard orbital representative,
$\theta$ generates $\co_F$ over $\co_v$.
Let $\phiw_0$ be the characteristic function of $W_{\co_v}$.
In this case, we could see that
$\gpX_{\ti z}(\Phi_0,s_1,s_2)=
(1-q_v^{-2s_1})^{-1}\gpX_{\ti z}'(\phiw_0,s_1+s_2,s_2)$
where
\begin{equation*}
\gpX_{\ti z}'(\phiw,a,b)
=	\int_{(\mk_v)^2\times k_v}
		|t_1|_v^a|t_2|_v^b
		\phiw\left(t_1,t_1{\cT}(u+t_2\theta),t_1{\cN}(u+t_2\theta)\right)
	\md t_1\md t_2 du
\end{equation*}
for $\phiw\in\sS(W_\A)$.
Let $\phiw_0'\in\sS(W_{\A})$ be the characteristic function
of $W_{\co_v}\setminus(\gp_v^2\times\gp_v\times\co_v)$.
Then the similar observation as in the case of (sp), we have
$\gpX_{\ti z}'(\phiw_0,a,b)=(1-q^{1-2a+b})^{-1}\gpX_{\ti z}'(\phiw_0',a,b)$.
Let $j=0$ if $\ti z$ is of type (ur) and $j=1$ if $\ti z$ is of type (rm).
Then it is easy to see that
$\phiw_0'(t_1,t_1{\cT}(u+t_2\theta),t_1{\cN}(u+t_2\theta))=1$
if and only if
\begin{equation*}
t_1\in\co_v^\times, t_2\in \co_v,u\in \co_v,
\qquad\text{or}\qquad
t_1\in\pi\co_v^\times, t_2\in\pi^{-j}\co_v,u\in\co_v.
\end{equation*}
Hence $\gpX_{\ti z}'(\phiw_0',a,b)=(1-q^{-b})^{-1}(1+q^{-a+jb})$
and this  finish the proof.
\end{proof}
As a corollary we obtain the following.
\begin{lem}\label{localconv}
For any $v\in\gM$, the local zeta function
$\gpX_{\ti y,v}(\phiu_v,s_1,s_2,\omega_{1v},\omega_{2v})$
is holomorphic in the region
$\Re(s_1)>0,\Re(s_2)>0,\Re(s_1+s_2)>1/2$.
\end{lem}
\begin{proof}
The statement for $v\in\gM_\infty$ follows from
Proposition \ref{prop:localfneq} (1).
For $v\in\gMf$, the support of $\phiu_v$ is contained in $\pi^{-m}U_{\co_v}$
for some integer $m$. Hence the result follows from Proposition
\ref{explicitformula} and the relatively invariant property
of the local zeta function under the action of $H_{k_v}$.
\end{proof}

\subsection{Orbital zeta functions}
We now discuss the relation between global and local situations.
\begin{defn}\label{defn:orbitalzeta}
For $\ti y\in U_{k}^\sst$, $\phiu\in \sS(U_\A)$
and $\omega_{1},\omega_{2}\in\Omega^1$,
we define
\begin{equation*}
X_{\ti y}^\ast(\phiu,s_1,s_2,\omega_{1},\omega_{2})
=
	\int_{H_\A}
	|t|^{s_1}|t\det g|^{2s_2}
	\omega_{1}(t)
	\omega_{2}(t\det g)
	\phiu(h\ti y)\, d_\pr h\\
\end{equation*}
and call it {\em the orbital zeta function}.
\end{defn}
Note that this integral depends only on $H_k$-orbits.
By the standard consideration
the global zeta function decompose into as follows
and this is the reason why we are interested in 
the orbital zeta functions.
\begin{lem}
We define
$X^\ast(\phiu,s_1,s_2,\omega_1,\omega_2)
=\Delta_k^{1/2}\gC_k^3\cdot
X(\phiu,s_1,s_2,\omega_1,\omega_2)$.
Then
\begin{equation*}
X^\ast(\phiu,s_1,s_2,\omega_{1},\omega_{2})
=2^{-1}\sum_{\ti y\in H_k\backslash U_k^\sst}
X_{\ti y}^\ast(\phiu,s_1,s_2,\omega_{1},\omega_{2}).
\end{equation*}
\end{lem}

The orbital zeta function has an Euler product.
We consider this Euler product more precisely.
For the rest of this subsection, we assume $\phiu$ is of the product form
$\phiu=\prod_{v\in\gM}\phiu_v$.
Let $\omega_i=\prod_{v\in\gM}\omega_{iv}$.
Then by definition
\begin{equation*}
X_{\ti y}^\ast(\phiu,s_1,s_2,\omega_{1},\omega_{2})
=\prod_{v\in\gM}\gpX_{\ti y,v}(\phiu_v,s_1,s_2,\omega_{1v},\omega_{2v}).
\end{equation*}
For each $v\in\gM$, let $\ti z_{\ti y,v}$ be the
standard representative lying in the orbit of $\ti y$
and put
\begin{equation*}
\Theta_{\ti y,v}(\phiu_v,s_1,s_2,\omega_{1v},\omega_{2v})
=\gpX_{\ti z_{\ti y,v},v}(\phiu_v,s_1,s_2,\omega_{1v},\omega_{2v}).
\end{equation*}
If we let $\ti y=h\ti z_{\ti y,v}$ for some $h=(t,g)\in H_{k_v}$,
then $t=R_1(\ti y)/R_1(\ti z_{\ti y,v})$
and $(t\det g)^2=R_2(\ti y)/R_2(\ti z_{\ti y,v})$.
We put
$\Delta_{\ti y,v}=R_2(\ti z_{\ti y,v})/R_2(\ti y)\in k_v^\times$.
Then we have
\begin{equation*}
\gpX_{\ti y,v}(\phiu_v,s_1,s_2,\omega_{1v},\omega_{2v})
=	\ti\omega_{1v}\left(\frac{R_1(\ti z_{\ti y,v})}{R_1(\ti y)}\right)
	\ti\omega_{2v}(\sqrt{\Delta_{\ti y,v}})
	\Theta_{\ti y,v}(\phiu_v,s_1,s_2,\omega_{1v},\omega_{2v}).
\end{equation*}
where we put 
$\ti\omega_{1v}=|\cdot|^{s_1}\omega_{1v}$
and $\ti\omega_{2v}=|\cdot|^{2s_2}\omega_{2v}$.
Note that Assumption \ref{omega2}
vanishes the ambiguity of the choice of the square root.
We put $\Delta_{\ti y}=(\Delta_{\ti y,v})_{v\in\gM}\in\ma$.
Since $\Delta_{\ti y}=(R_2(\ti z_{\ti y,v}))_{v\in\gM}/R_2(\ti y)$,
by the observation after Definition \ref{standardrepresentative},
if we regard $\Delta_{\ti y}$ and $\Delta_{\ti k(\ti y)/k}$
as elements of $\ma/\mk$ then they coincide.
Consider the product of the above formula for $v\in\gM$.
Since $R_1(\ti z_{\ti y,v})=1$ for all $v\in\gM$ and $\omega_1(R_1(\ti y))=1$,
we have
\begin{equation*}
X_{\ti y}^\ast(\phiu_v,s_1,s_2,\omega_{1v},\omega_{2v})
=	\ti\omega_{2}(\sqrt{\Delta_{\ti y}})
	\prod_{v\in\gM}\Theta_{\ti y,v}(\phiu_v,s_1,s_2,\omega_{1v},\omega_{2v}).
\end{equation*}
For a finite set $T$ of places of $k$,
we define the truncated zeta function
$\zeta_{k,T}(s)=\prod_{v\in \gMf\setminus T}(1-q_v^{-s})^{-1}$.
For an separable quadratic algebra $F$ of $k$,
we define $\zeta_{F,T}(s)$ similarly.
By Lemmata \ref{localcont}, \ref{localconv}
and Proposition \ref{explicitformula},
we have the meromorphic continuation of the orbital zeta functions
as below.
\begin{prop}\label{prop:orbitalzeta}
Let $T\supset\gM_\infty$ be a finite set
such that $\phiu_v$ is the characteristic function of $U_{\co_v}$
and $\omega_v$ is unramified unless $v\in T$.
Then
\begin{multline*}
X_{\ti y}(\phiu,s_1,s_2,\omega_1,\omega_2)
=	N(\Delta_{\ti k(\ti y)/k})^{-s_2}
	\omega_2(\sqrt{\Delta_{\ti y}})
	\prod_{v\in T}
		\Theta_{\ti y,v}(\phiu_v,s_1,s_2,\omega_{1v},\omega_{2v})\\
\times
	\zeta_{k,T}(2s_2)\zeta_{k,T}(2s_1+2s_2-1)\zeta_{k,T}(s_1+2s_2)
	\zeta_{k(\ti y),T}(s_1)\zeta_{k(\ti y),T}(s_1+2s_2)^{-1}.
\end{multline*}
This function is meromorphic on $\C^2$
and holomorphic in the region $\Re(s_1)>1,\Re(s_2)>1$.
\end{prop}

\end{document}